\documentclass[12pt, onecolumn]{IEEEtran}
\usepackage{graphicx}
\usepackage{latexsym}
\usepackage{epsfig}
\usepackage{grffile}
\usepackage{balance}
\usepackage{cite}
\usepackage{array}
\usepackage{mdwtab}
\usepackage{eqparbox}
\usepackage{units}
\usepackage{amsmath}
\usepackage{color}
\usepackage{amssymb}
\usepackage{float}
\usepackage{caption}
\usepackage{setspace}
\usepackage{algorithmicx}
\usepackage{algpseudocode}
\usepackage{epstopdf}
\usepackage{soul}
\usepackage[font=footnotesize]{caption}
\begin{document}

\title{Cooperative  Beamforming  for  Cognitive Radio-Based Broadcasting Systems with Asynchronous Interferences}
\author{%
\IEEEauthorblockN{Mai H. Hassan\IEEEauthorrefmark{2} and Md. Jahangir Hossain\IEEEauthorrefmark{4}\\}
\IEEEauthorblockA{
\IEEEauthorrefmark{2}Department of Electrical and Computer Engineering\\The University of British Columbia, Vancouver, BC, Canada\\ 
\IEEEauthorrefmark{4}School of Engineering\\ The University of British Columbia, Kelowna, BC, Canada\\
\emph{maih@ece.ubc.ca}, \emph{jahangir.hossain@ubc.ca}
}
\thanks{This work is supported by the Natural Science and Engineering Research Council (NSERC) of Canada, and will be presented in part at the IEEE Wireless Communications and Networking Conference (WCNC'14), Istanbul, Turkey, April 2014.}
}%

\maketitle

\begin{abstract}
In a cooperative cognitive radio (CR) network, 
cooperative beamforming 
 can enable  concurrent transmissions of both primary  and secondary systems at a given channel. 
However, such cooperative beamforming can introduce  asynchronous interferences at the primary receivers as well as at the secondary receivers and these asynchronous interferences are overlooked in beamforming design. 
In order to address the asynchronous interference issue for a generalized scenario with multiple primary and multiple secondary receivers, in this paper, we  propose an innovative cooperative beamforming technique. In particular,  the cooperative beamforming design is formulated as an optimization problem that maximizes the weighted sum achievable transmission rate  of secondary destinations while it  maintains the asynchronous interferences at the primary receivers below their target thresholds. In light of the intractability of the
problem,  
we propose a two-phase suboptimal  cooperative beamforming technique.  First, it  finds the beamforming directions corresponding to different secondary destinations.   Second,  it allocates the power among different beamforming directions.  Due to the multiple interference constraints corresponding to multiple primary receivers,  the power allocation scheme in the second phase is still  complex. Therefore, we also propose a low complex power allocation algorithm.  The proposed beamforming technique is  extended for the cases, when  cooperating CR nodes (CCRNs) have statistical or erroneous channel knowledge of  the primary receivers. We also investigate the performance of joint CCRN selection and  beamforming technique.
 The presented numerical results show that the proposed beamforming technique can significantly reduce the asynchronous interference signals at the  primary receivers and increase the sum transmission rate of secondary destinations compared to  the well known  zero-forcing beamforming (ZFBF) technique.  
%

\end{abstract}


\section{Introduction}

Recently  dynamic spectrum access (DSA) or opportunistic spectrum access policy has received a great deal of attention in order to improve the overall spectrum utilization. Cognitive radio  (CR) \cite{mitola},  \cite{haykin}  is one of the key enabling technologies  in order to facilitate DSA.  
 Different DSA  mechanisms  have already   been envisioned and studied  in the literature  \cite{Ekram_book},  \cite{DSA_Survey}.  Among these, two  approaches namely,  underlay  and  overlay spectrum access mechanisms for spectrum sharing between primary and secondary systems have been   considered  widely. 
 The underlay spectrum access mechanism \cite{DSA_Survey}   allows simultaneous sharing of underutilized frequency
bands  by a secondary/CR system along with a primary system provided that the introduced interference to the primary users does not exceed certain thresholds specified  by the primary system or the regulatory authority,  see for example  \cite{interference_threshold} for details. 

Cooperation among nodes   in a wireless network can improve the overall  performances \cite{cooperative_networks1,cooperative_networks2}.   For example,   small nodes with simple omni-directional antenna can cooperatively emulate a large highly directional antenna array which is referred to as cooperative transmit beamforming.  In other words, in order to send a common message, a number of single antenna-based  nodes in a wireless network organize themselves into a virtual antenna array and focus their transmission in the direction of the intended receivers.  Such beamforming technique has been proposed and studied for traditional wireless communication networks  e.g.,  wireless sensor network as this potentially offers large increases in energy efficiency, in attainable range and transmission rate, see for example \cite{dist_beamforming}  and the references therein.     With such cooperative beamforming technique,  the achievable  data rate gain is quite compelling in spite of  certain   costs  associated with it  e.g.,  synchronising the sensor nodes and  the local exchange of sensor nodes' observations \cite{Madhow_Dist_BF}.

  
   In order to reap the benefit of cooperative communications, cooperative beamforming  technique has been  proposed for CR systems as well \cite{distributed_beamforming, ref1_coop, ref2_coop}. 
 In fact  cooperative transmit beamforming  can be very effective for  CR systems that work based on underlay spectrum access  mechanism which  imposes severe constraints on the transmission power of CR systems \cite{DSA_Survey}. For example, if a cognitive/secondary transmitter\footnote{In the literature CR users are also referred to as secondary users and throughout this paper the words secondary and cognitive have been used interchangeably.}  wants to broadcast a common message to a group  of secondary destinations,  the transmitter  may not be allowed to transmit  enough power 
  to cover all the destinations  due to the interference restriction  imposed by the nearby primary system.       In such situation,    a group of  secondary nodes which are referred to as cooperating cognitive radio nodes (CCRNs)  can collaboratively use transmit beamforming to broadcast the common message to the secondary destinations see for example \cite{distributed_beamforming}.  
%
   Using an innovative orthogonal projection technique, the authors in \cite{distributed_beamforming} obtained the so-called zero forcing beamforming (ZFBF) weights of the CCRNs  to null the interference at the primary receivers. Using the ZFBF technique, the authors in \cite{ref1_coop} proposed a cross-layer optimization of the transmission rate and scheduling scheme of the data packets at the secondary source and at the CCRNs in the CR network. In \cite{ref2_coop}, power allocation for cooperative CR networks was studied, along with user selection under imperfect spectrum sensing.


       Given the fact that in practice, different CCRNs  are usually located in different geographical locations, their signals can arrive with different propagation delays at each primary receiver and at each  secondary receiver.  Therefore with such cooperative beamforming, simultaneous transmissions from the CCRNs can cause {\it  asynchronous interferences} which are discussed in details in Section \ref{asynch}. 
    Although the cooperative beamforming technique can improve the overall performance of CR systems,  the  asynchronous  interferences  are  overlooked in designing beamforming technique.   As we will see later in this paper that  the  ZFBF  technique introduces a significant amount of asynchronous interference power at the primary receivers. 
   
    The asynchronous interference issue for conventional cooperative multi-cell mobile networks  has been studied in \cite{asynch_mitigation}, where multiple base stations (BSs)  cooperate together to simultaneously transmit information to  each mobile user in the network. 
 However in cooperative beamforming for CR systems, the design goal is different and  a new beamforming technique is required.  In our earlier  work \cite{ my_journal} (the conference version published in \cite{my_icc})   considering the  asynchronous interference, we developed a cooperative beamforming technique for a CR network with \emph{one} secondary destination and \emph{one} primary receiver. 

 As a follow-up of our initial  work \cite{my_journal},   in this paper we consider a more generalized setup where a group of CCRNs uses cooperative beamforming technique to broadcast common message to \emph{multiple}  secondary destinations using a wireless communication channel that is used by a primary transmitter  to transmit information to \emph{multiple} primary receivers simultaneously. It is also considered that these primary receivers have different interference constraints  in general. 
With multiple primary and secondary receivers, asynchronous interferences are   introduced not only at the primary receivers but also at the secondary destinations except one secondary destination.  Due to   these asynchronous interferences, the optimal beamforming technique developed in  our earlier work \cite{my_journal}  can not be extended for a generalized system with multiple primary and multiple secondary receivers.   In fact, the optimal beamforming techique is intractable due to the non-convexity and non-linearity of the problem. Even then development of suboptimal beamforming technique is complex due to multiple interference constraints corresponding to multiple primary receivers  which is discussed later.    In order to address the asynchronous interference issues for such a generalized scenario,  in this paper,  we develop  innovative cooperative beamforming techniques.  In particular, the contributions of this paper can be summarized as follows:  
\begin{itemize}
\item For a generalized scenario,  the cooperative beamforming design is formulated as an optimization problem that maximizes the weighted sum achievable transmission rate  of secondary destinations while it maintains the interference thresholds at the primary receivers. Due to the non-convexity and non-linearity  of formulated optimization problem, we propose a two-phase suboptimal beamforming technique. First, it finds the beamforming direction corresponding to a  secondary destination  that maximizes the received signal power at that  secondary destination while it minimizes  the asynchronous interference power at other secondary destinations and at all primary receivers.   Second, it allocates the power among different beamforming directions to maintain the interference constraints at the primary receivers. Due to the multiple interference constraints, the power allocation scheme in second phase can be  complex as discussed later. Therefore, we also propose a low complex power allocation  (LCPA) scheme.  The  presented numerical results show that the developed cooperative beamforming technique can increase the sum  data rate of the secondary destinations  up to 64\% compared to the well-known ZFBF technique.  
\item 
%
 We extend the proposed cooperative beamforming technique to the case of having only partial channel state information (CSI)   between the primary receivers and the CCRNs. This partial CSI has been modeled by two scenarios. In the first scenario, we consider having only statistical CSI  of the channels. In this case, the asynchronous interferences at the primary receivers are guaranteed in a statistical sense \cite{kim}, \cite{stat_channel}.  In the absence of mathematically tractable expression of   the distribution  of the random interference power  at the primary receiver,   we develop an upper bound on  the  probability of introducing interference at a given primary receiver beyond a given threshold.   Then this developed upper bound is used to design a robust leakage beamforming (RLBF) technique.    The second scenario considers of  having erroneous CSI. For both scenarios, we design  RLBF techniques that can protect the primary network's  functionality by satisfying the interference constraints at all primary receivers in the network, in spite of the partial  knowledge of CSI. 
\item We also investigate the  performance of joint CCRN selection and  beamforming technique. The numerical results show that CCRN selection in conjunction with beamforming can further increase the sum transmission rate of secondary destinations significantly. 
\end{itemize}


The rest of the paper is organized as follows. In Section \ref{sec:system_model} we present the overall system description and model the asynchronous interference signals at the primary receivers as well as at the secondary destinations  mathematically. While  in Section III  we develop the beamforming techniques with perfect CSI at the CCRNs,  in Section IV,  we develop robust beamforming  techniques for  two scenarios, i.e., imperfect channel CSI and statistical CSI. In  Section V, we investigate the performance of joint CCRN selection and cooperative beamforming.  Section \ref{sec:Numerical-results}  presents some numerical examples.  Finally, Section \ref{sec:conc} concludes the paper.

\section{\label{sec:system_model}System Model }
\subsection{Overall Description and Operating Assumptions}
\begin{figure}[h]
\centering
\includegraphics[scale=0.5,angle=-90,trim=2cm 1cm 2.5cm 0cm]{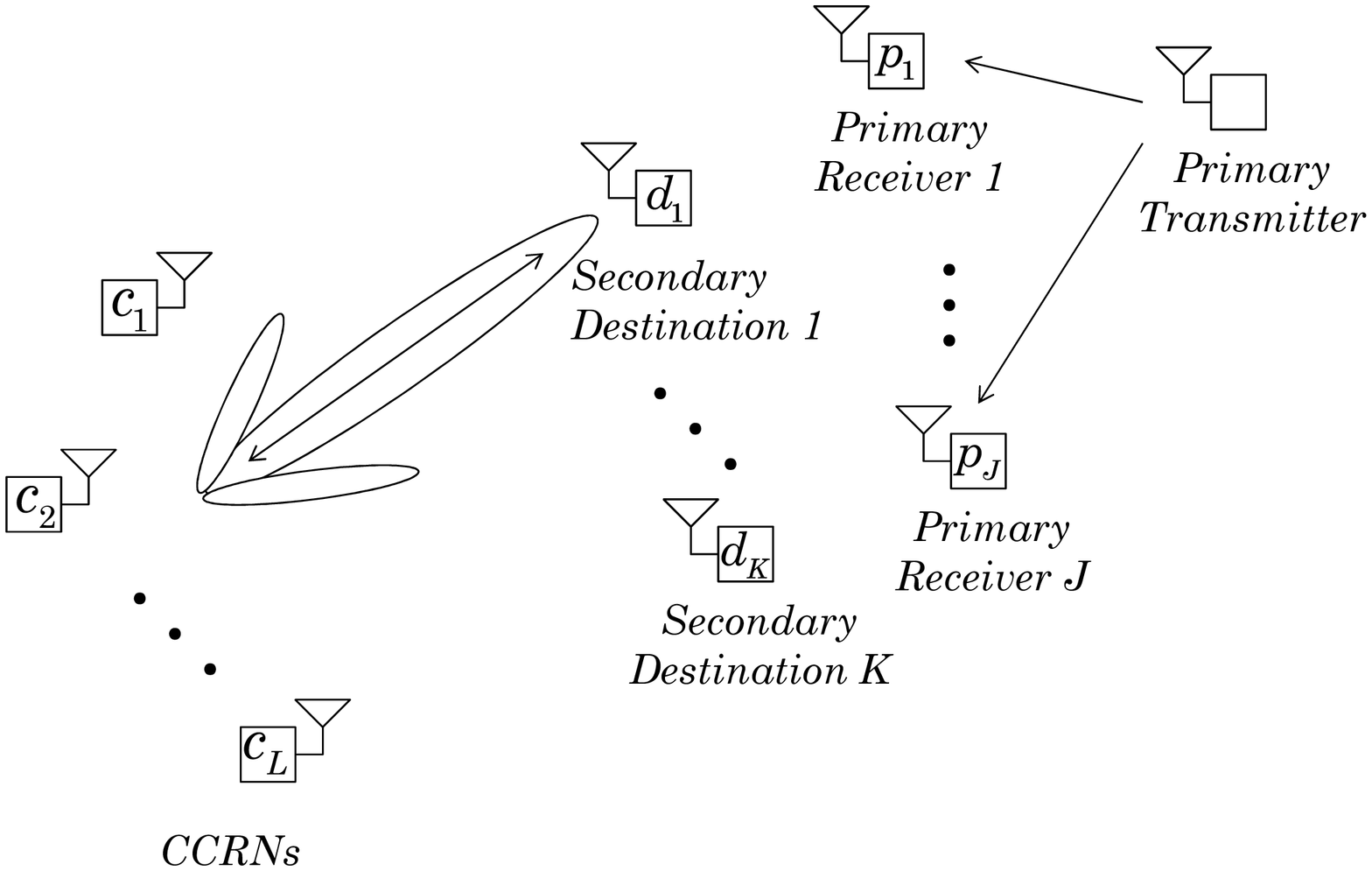}
\caption{System model for cooperative beamforming with $L$ CCRNs, $K$ secondary destinations and $J$ primary receivers.}
\label{system_model}
\end{figure}

For an example, the cooperative beamforming for CR-based broadcasting system is shown in  Fig. \ref{system_model}  where  a group of $L$ CCRNs  uses a transmit beamforming technique to transmit  common information to a group of $K$   secondary destinations.  All nodes  are assumed to be equipped with single antenna.  As we mentioned earlier that   similar type of cooperative beamforming scenario is considered for traditional wireless networks e.g., wireless  sensor network  due to its compelling gain in the transmission rate, see for example \cite{dist_beamforming}, \cite{Madhow_Dist_BF}.   
Using the underlay spectrum access mechanism, the CR system   shares a communication channel with   a primary transmitter e.g., a primary BS   that transmits  information to $J$ primary receivers simultaneously.  
For notational convenience,  $K$ secondary destinations are denoted by d$_k, k=1,\cdots,K$, $L$ CCRNs are denoted by c$_l,~l=1,\cdots,L$ and $J$ primary receivers are denoted by p$_{j},~j=1,\cdots,J$. 
In what follows, we provide the operating principles as well as the assumptions that we consider in our problem formulation.  



\begin{figure}[h]
\centering
\includegraphics[scale=0.8,trim=0cm 18.5cm 0cm 2cm]{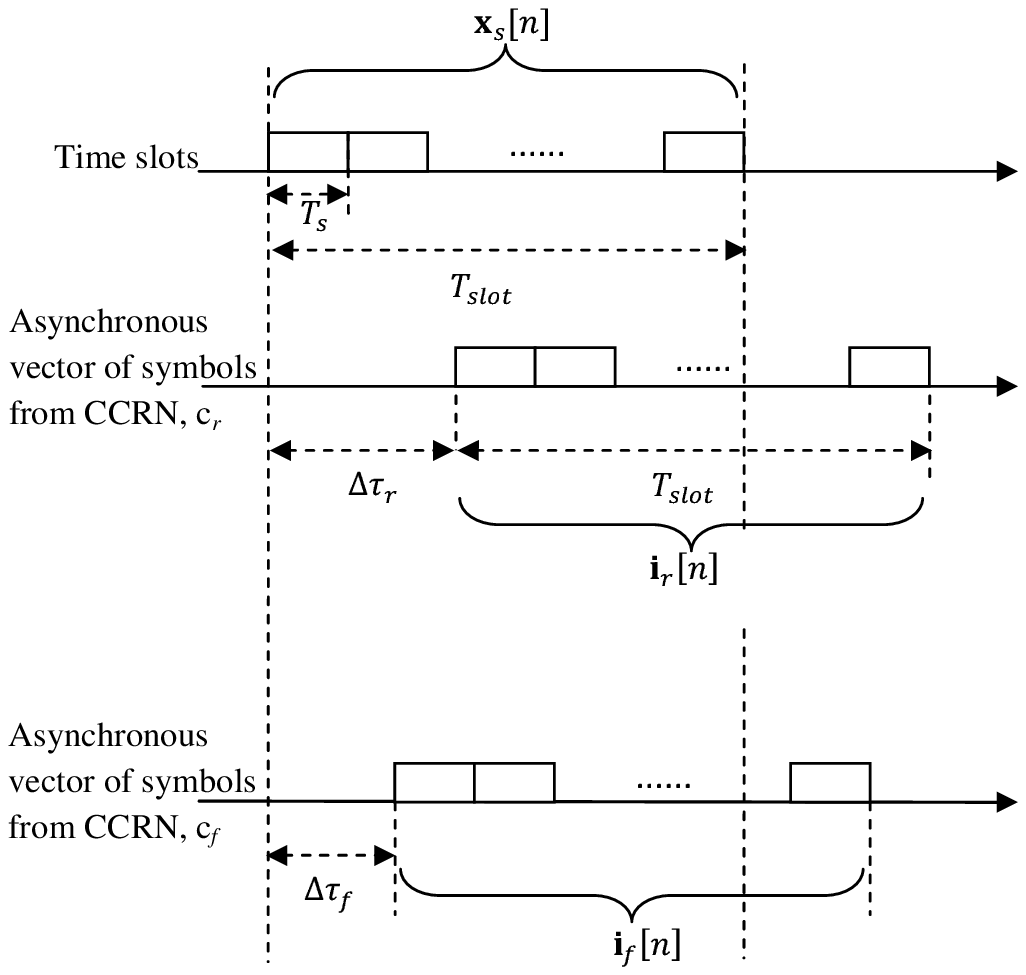}
\caption{An example of the asynchronous vector of symbols received at each of the primary receivers, as well as other secondary destinations. $T_s$ is the symbol duration and $T_\text{slot}$ is the time slot duration.\label{fig:The-asynchronous-interference}}
\end{figure}

We consider that  both primary and secondary systems work in a time-slotted fashion with slot duration $T$ sec.  
At CCRN, c$_l$  the information stream is  mapped into modulated symbols, $x_{s}$ which has average power $P$  and  the data vector consisting of these $M$ modulated symbols is denoted by $\mathbf{x}_{s}$. We assume a block fading channel model, in which the channel fading is assumed to remain roughly the same over the time slot, but is independent of the fading in other time slots.
 The set of cooperative CCRNs  uses $K$ different beamforming vectors to direct transmission  to  $K$ different secondary destinations. 
The received signal at secondary destination, d$_k$   can be written as 
\begin{equation}
\mathbf{y}_{k}[n]=\mathbf{h}_{k}^{\mbox{\tiny{s}}}[n]\mathbf{g}_{k}[n]\mathbf{x}_s[n]+ \mathbf{I}_k[n]+\mathbf{m}_k[n]+\mathbf{z}_{k}[n],\end{equation}
where $\mathbf{x}_{s}[n]$ is common message symbols transmitted at time slot $n$ and  $\mathbf{h}_{k}^{\mbox{\tiny{s}}}[n]\triangleq [h_{k1}^{\mbox{\tiny{s}}}[n],\text{\ldots},h_{kL}^{\mbox{\tiny{s}}}[n]]$ is the channel vector from $L$  transmitting CCRNs  to secondary destination, d$_k$.  The vector  $\mathbf{g}_{k}[n]\triangleq [g_{k1}[n],\text{\ldots},g_{kL}[n]]^{T}$ denotes  the beamforming weight vector of the set of CCRNs corresponding  to transmission  to  secondary destination, d$_k$  with each element $g_{kr}$ denoting the weight of the CCRN, c$_r$.  $\mathbf{z}_{k}[n]$ is the additive white Gaussian noise (AWGN) vector at  secondary destination, d$_k$  with zero mean and two sided power  spectrum density $N_0/2$,  and $\mathbf{m}_k[n]$  is the received interfering signal vector from the primary BS at secondary destination d$_k$. $\mathbf{I}_k[n]$ is  the  asynchronous interference signal at  secondary destination, d$_k$  resulting from the data transmissions to the other $(K-1)$ secondary destinations and can be written  as follows
\begin{equation}
\mathbf{I}_k[n]=\sum_{i=1, i\neq k}^{K}\sum_{r=1}^{L}h_{kr}^{\mbox{\tiny{s}}}[n]g_{ir}[n]\mathbf{i}_{r}[n],
\end{equation}
where $\mathbf{i}_{r}[n]$ is the asynchronous vector of symbols received at the secondary destination d$_k$  from the CCRN  c$_r$, as shown in Fig. \ref{fig:The-asynchronous-interference}.

The ZFBF technique developed in \cite{distributed_beamforming}  did not consider the asynchronous interference issue described below. In fact, we will see later in this paper  that  if ZFBF technique is used in this scenario, the asynchronous interferences at the primary receivers exceed the target thresholds.  So there is a need for developing an innovative beamforming technique which is the main focus of this paper. 

\subsection{\label{asynch}Modeling of Asynchronous Interferences}
Given the fact that the CCRNs  are  located in different geographical locations, the received signals from different CCRNs  at different primary receivers and at different secondary destinations can experience different propagation delays.  
Although the received signal at a particular secondary destination e.g.,  d$_1$ from different CCRNs  can be synchronized by using a timing advance mechanism which is currently employed  in the uplink of GSM and 3G cellular networks (see for example \cite{asynch_mitigation}) or other mechanism  \cite{Madhow_Dist_BF}, the received signals at the primary receivers p$_j$ ($j=1,\cdots,J$)  and at the other secondary destinations, d$_k$ $(k=2,\cdots, K)$  cannot be synchronized simultaneously.  
As such, the signal transmissions from CCRNs will introduce   asynchronous interference at  primary receivers p$_j$ ($j=1,\cdots,J$)  and at the other secondary destinations, d$_k$ $(k=2,\cdots, K)$. In our previous work  \cite{my_journal}, we have modeled the asynchronous interference  at the primary destination with one primary receiver and one secondary destination in the system. However  for the generalized scenario considered in this paper,  we need to model   the asynchronous interferences not only at different primary receivers but also at different secondary destinations. In what follows, we model these asynchronous interferences. For notational convenience, we will drop the time slot index $n$.

\subsubsection{Asynchronous interference at  primary receiver, p$_j$}
Mathematically,   the asynchronous interference power resulting from transmission to  secondary destination, d$_k$  at  primary receiver, p$_j$,   $P_\text{asynch}^{(j,k)}$, can be written in the following  form \cite{my_journal}
\begin{equation}
P_\text{asynch}^{(j,k)}=\sum_{r=1}^{L}\sum_{f=1}^{L}g_{kf}^{\dagger}({h_{jf}^{\mbox{\tiny{p}}}})^{\dagger}h_{jr}^{\mbox{\tiny{p}}} g_{kr}{\beta_k^{j}}^{(r,f)},
\label{AsynPower}
\end{equation}
where $h_{jr}^{\mbox{\tiny{p}}}$ is the  channel fading gain  from CCRN c$_r$ to  primary receiver, p$_j$,   ${\beta_k^{j}}^{(r,f)}$ is the correlation between the asynchronous symbols of  CCRNs, c$_r$ and c$_f$  at  primary receiver, p$_j$  corresponding to  the transmission to  secondary destination, d$_k$.  The value of ${\beta_k^{j}}^{(r,f)}$ can be calculated for given  propagation delays between CCRNs, c$_r$ and c$_f$ to  primary receiver, p$_j$  using the same technique described in  \cite{my_journal}.
Now the total  asynchronous interference power at   primary receiver, p$_j$  can be expressed as  
\begin{equation}
P^{j}_\text{asynch}=\sum_{k=1}^{K}\sum_{r=1}^{L}\sum_{f=1}^{L}g_{kf}^{\dagger}(h_{jf}^{\mbox{\tiny{p}}})^{\dagger}h_{jr}^{\mbox{\tiny{p}}}g_{kr}{\beta_k^{j}}^{(r,f)}.
\label{eq:leak_power}\end{equation}
The asynchronous interference power at primary receiver, p$_j$  in eq. (\ref{eq:leak_power}) can be rewritten in a matrix form as follows
\begin{equation}
P^{j}_\text{asynch}=\sum_{k=1}^{K}{\mathbf{g}_k}^{\dagger}\mathbf{R}_k^{j}\mathbf{g}_k\label{eq:leak_power2}\end{equation}
where $\mathbf{R}_k^{j}$ is 
 expressed as 
\begin{equation}
\label{chCovMatrix}
\mathbf{R}_k^{j}=
\left[\begin{array}{ccc}
{\beta_k^{j}}^{(1,1)}(h_{j1}^{\mbox{\tiny{p}}})^{\dagger}h_{j1}^{\mbox{\tiny{p}}} & \negthickspace\cdots\negthickspace & {\beta_k^{j}}^{(1,L)}(h_{j1}^{\mbox{\tiny{p}}})^{\dagger}h_{jL}^{\mbox{\tiny{p}}}\\
{\beta_k^j}^{(2,1)}(h_{j2}^{\mbox{\tiny{p}}})^{\dagger}h_{j1}^{\mbox{\tiny{p}}} & \negthickspace\cdots\negthickspace & {\beta_k^j}^{(2,L)}(h_{j2}^{\mbox{\tiny{p}}})^{\dagger}h_{jL}^{\mbox{\tiny{p}}}\\
\negthickspace\vdots\negthickspace & \negthickspace\ddots\negthickspace & \negthickspace\vdots\negthickspace\\
{\beta_k^j}^{(L,1)}(h_{jL}^{\mbox{\tiny{p}}})^{\dagger}h_{j1}^{\mbox{\tiny{p}}} & \negthickspace\cdots\negthickspace & {\beta_k^j}^{(L,L)}(h_{jL}^{\mbox{\tiny{p}}})^{\dagger}h_{jL}^{\mbox{\tiny{p}}}\end{array}\right].  \end{equation}

The received signal power at  secondary destination, d$_k$  is given by
\begin{equation}
P_{k,\text{signal}}=P{\mathbf{g}_k}^{\dagger}(\mathbf{h}_k^{\mbox{\tiny{s}}})^{\dagger}\mathbf{h}_k^{\mbox{\tiny{s}}}\mathbf{g}_k. \label{eq:signal_power}\end{equation}

\subsubsection{Asynchronous interference at  secondary destination, d$_k$}

The asynchronous interference power at  secondary destination, d$_k$  resulting from  transmission  to  secondary destination, d$_i$  is given by
\begin{equation}
  \mbox{AI}^{k}_i=\sum_{r=1}^{L}\sum_{f=1}^{L}{\beta^{k}_{i}}^{(r,f)}g_{if}^{\dagger}(h_{kf}^{\mbox{\tiny{s}}})^{\dagger}h_{kr}^{\mbox{\tiny{s}}}g_{ir},
\end{equation}
where ${\beta_i^{k}}^{(r,f)}$ is the correlation between the asynchronous symbols of  CCRNs, c$_r$ and c$_f$  at  secondary destination, d$_i$  corresponding to  the transmission to  secondary destination, d$_k$ where $k\neq i$.  The value of ${\beta_i^{k}}^{(r,f)}$ can be calculated for given  propagation delays between CCRNs, c$_r$ and c$_f$ to  secondary destinations, d$_i$ and d$_k$  using the same method described in  \cite{my_journal}.  
Therefore, the total asynchronous interference power at  secondary destination, d$_k$ resulting from  data transmission to the other $(K-1)$ secondary destinations,  {AI}$^k$ can be written as 
\begin{equation}
 \mbox{AI}^k=\sum_{
i=1, i\neq k}^{K}\sum_{r=1}^{L}\sum_{f=1}^{L}{\beta^{k}_{i}}^{(r,f)}g_{if}^{\dagger}(h_{kf}^{\mbox{\tiny{s}}})^{\dagger}h_{kr}^{\mbox{\tiny{s}}}g_{ir}.
\end{equation}
Similar to eq. (\ref{eq:signal_power}),   {AI}$^k$  can be written in a matrix form as follows 
\begin{equation}
\mbox{AI}^k=\sum_{
i=1, i\neq k}^{K}\mathbf{g}_{i}^{\dagger}\mathbf{T}^k_{i}\mathbf{g}_{i},
\end{equation}
where $\mathbf{T}^k_{i}$ is  written as  
\begin{equation}
\mathbf{T}^k_{i}\triangleq\left[\begin{array}{ccc}
{\beta^{k}_{i}}^{(1,1)}(h_{k1}^{\mbox{\tiny{s}}})^{\dagger}h_{k1}^{\mbox{\tiny{s}}} & \cdots & {\beta^{k}_{i}}^{(1,L)}(h_{k1}^{\mbox{\tiny{s}}})^{\dagger}h_{kL}^{\mbox{\tiny{s}}}\\
\vdots & \ddots & \vdots\\
{\beta^{k}_{i}}^{(L,1)}(h_{kL}^{\mbox{\tiny{s}}})^{\dagger}h_{k1}^{\mbox{\tiny{s}}} & \cdots & {\beta^{k}_{i}}^{(L,L)}(h_{kL}^{\mbox{\tiny{s}}})^{\dagger}h_{kL}^{\mbox{\tiny{s}}}\end{array}\right].
\nonumber
\end{equation}
 
\section{\label{perfect_channel} Beamforming Design with Perfect Channel Knowledge}

In this section, we develop a new beamforming technique, called cooperative leakage beamforming (LBF)  technique in order to address the problem of asynchronous interferences at the primary receivers and other secondary destinations. In this development, we use the same assumption as in \cite{distributed_beamforming, ref3_globecom2010, ref5_globecom2010, ref6_globecom2010} that the channel fading gains, i.e., instantaneous CSI  between the CCRNs  and the primary receivers as well as the instantaneous CSI  between the CCRNs  and the secondary destinations are  known perfectly at the CCRNs.  Different possible scenarios have been considered in the literature in order to estimate the CSI  between the CCRNs and the primary receivers (see for examples, \cite{cooperation_with_PU},  \cite{Lampe_paper}).
 In the next section,  we consider the case when the CSI  between the CCRNs and the primary receivers are not known perfectly. 
\subsection{Problem Formulation}
The achievable  transmission rate of  secondary destination d$_k$,   $r_{k}$  can be expressed using the ideal capacity formula  as follows
\begin{equation}
r_{k}=\log_2\Biggl(1+\frac{P\mathbf{g}_{k}^{\dagger}(\mathbf{h}_{k}^{\mbox{\tiny{s}}})^{\dagger}\mathbf{h}_{k}^{\mbox{\tiny{s}}}\mathbf{g}_{k}}{\sigma_{k}^{2}+{\displaystyle \sum_{
i=1, i\neq k}^{K}\mathbf{g}_{i}^{\dagger}\mathbf{T}^k_{i}\mathbf{g}_{i}}}\Biggr),\end{equation}
where $\sigma_{k}^{2}$ is the power of the AWGN   plus the interference power from the primary transmitter at  secondary destination d$_k$. 
%
The   goal  is to design $K$ different beamforming vectors corresponding to $K$ secondary destinations that  maximize the weighted sum rate of all secondary destinations  while keeping the interference to the primary receivers below their target thresholds. We consider maximizing the weighted sum rate of all the $K$ secondary destinations since it is more  generalized (see for example \cite{weighted},  and the references therein).
The design goal can be formulated as an optimization problem as follows\footnote{We do not consider transmit power constraint for the CCRNs as we develop the beamforming technique for the interference limited scenario.}
\begin{eqnarray}
\mathbf{g}_{1}^\text{\tiny{(opt)}},\mathbf{g}_{2}^\text{\tiny{(opt)}},\cdots,\mathbf{g}_{K}^\text{\tiny{(opt)}} & = &  \max_{\mathbf{g}_{1},\cdots,\mathbf{g}_{K}}\sum_{k=1}^{K}w^k \log_2\Biggl(1+\frac{P\mathbf{g}_{k}^{\dagger}(\mathbf{h}_{k}^{\mbox{\tiny{s}}})^{\dagger}\mathbf{h}_{k}^{\mbox{\tiny{s}}}\mathbf{g}_{k}}{\sigma_{k}^{2}+{\displaystyle \sum_{
i=1, i\neq k}^{K}\mathbf{g}_{i}^{\dagger}\mathbf{T}^k_{i}\mathbf{g}_{i}}}\Biggr),\nonumber \\
\text{subject to: } & \text{ } & \sum_{k=1}^{K}\mathbf{g}_{k}^{\dagger}\mathbf{R}_k^j\mathbf{g}_{k}\leq\gamma^j_\text{th}, \text{for } j=1,\cdots, J.\label{optimum_form}\end{eqnarray}
where $w^k$ is the weighting factor of  secondary destination d$_k$,   and $\gamma^j_\text{th}$ is the required interference threshold for primary receiver, p$_j$. 


\subsection{Development of Suboptimal Cooperative LBF Technique}
\label{coopLBF}
The optimization problem in eq. (\ref{optimum_form}) is a non-linear and non-convex optimization problem due to the presence of the interference power AI$^k= \sum_{
i=1, i\neq k}^{K}\mathbf{g}_{i}^{\dagger}\mathbf{T}^k_{i}\mathbf{g}_{i}$   in  secondary destination d$_k$'s  transmission rate, $r_k$.   In light of the intractability of this optimization problem,  we propose a two-phase  suboptimal  cooperative LBF technique as described below. 

\subsubsection{\label{phase_I}  Phase  I} In this phase, we find  the direction of the normalized beamforming vector, $\mathbf{\bar{g}}_{k}$ that   maximizes the received signal power at   secondary destination d$_k$  while it minimizes  the interference at all primary receivers and  other secondary destinations. This  can be written as  the following optimization problem 
\begin{eqnarray}
\mathbf{\bar{g}}_{k}^\text{\tiny{(LBF)}} & = &    \max_{\mathbf{\bar{g}}_{k}} \frac{\mathbf{\bar{g}}_{k}^{\dagger}(\mathbf{h}_{k}^{\mbox{\tiny{s}}})^{\dagger}\mathbf{h}_{k}^{\mbox{\tiny{s}}}\mathbf{\bar{g}}_{k}}{\mathbf{\bar{g}}_{k}^{\dagger}\bigl(\mathbf{R}_k+\mathbf{T}_{k}\bigr)\mathbf{\bar{g}}_{k}}, \text{\;\;\;\;\;\;for } k=1,\cdots, K, \label{opt_dir}\end{eqnarray}
where $\mathbf{T}_{k}=\sum_{i=1, i\neq k}^{K} \mathbf{T}^i_{k}$  
and $\mathbf{R}_{k}=\sum_{j=1}^{J} \mathbf{R}^j_{k}$. 
The signal-to-leakage power ratio in eq. (\ref{opt_dir}) is in the form of a generalized Rayleigh quotient, that is maximized when $\mathbf{\bar{g}}_{k}^\text{\tiny{(LBF)}}$ is the normalized eigen vector of the matrix ${\bigl(\mathbf{R}_k+\mathbf{T}_{k}\bigr)}^{-1}(\mathbf{h}_{k}^{\mbox{\tiny{s}}})^{\dagger}\mathbf{h}_{k}^{\mbox{\tiny{s}}}$
that corresponds to its maximum eigen value  \cite{rayleigh_book}. As indicated before, the optimization problem in eq. (\ref{optimum_form}) is a non-linear and non-convex optimization problem which cannot be solved optimally, due to the presence of the asynchronous interference power {AI}$^k$.  By minimizing such interference and for mathematical tractability,   we neglect  the interference power at secondary destination d$_k$.  In Section \ref{sec:Numerical-results}, we will  show  that after minimizing the  asynchronous interference power at the secondary destinations, the remaining asynchronous interference power has a negligible effect on  transmission rate $r_k$. So by neglecting the asynchronous interference power  {AI}$^k$, the transmission rate $r_k$  can be now approximated as
\begin{equation}\label{rApprox}
r_{k}^{\mbox{\tiny{App}}}\approx\log_2\Bigl(1+\frac{P\mathbf{\alpha}_{k}\mathbf{\bar{g}}_{k}^{\text{\tiny{(LBF)}}{\dagger}}(\mathbf{h}_{k}^{\mbox{\tiny{s}}})^{\dagger}\mathbf{h}_{k}^{\mbox{\tiny{s}}})\mathbf{\bar{g}}_{k}^\text{\tiny{(LBF)}}}{\sigma_{k}^{2}}\Bigr),\end{equation}

where $\mathbf{\alpha}_k$ is the power allocated to the beamforming direction   corresponding to  secondary destination d$_k$.
\subsubsection{\label{multiple_SDs}Phase II}  In this phase, we propose to  allocate power $\mathbf{\alpha}_k^\text{\tiny{(LBF)}}$ among different beamforming directions.  As such the approximated weighted sum rate of secondary destinations is maximized while the interference thresholds at different primary receivers are met.  In particular, given the normalized  beamforming vector $\mathbf{\bar{g}}_{k}^\text{\tiny{(LBF)}}$ obtained in Phase-I,  we obtain its allocated power $\mathbf{\alpha}_k^\text{\tiny{(LBF)}}$ that  satisfies the interference threshold at all primary receivers simultaneously, where $\mathbf{{g}}_{k}^\text{\tiny{(LBF)}}=\sqrt{\alpha_k^\text{\tiny{(LBF)}}}\mathbf{\bar{g}}_{k}^\text{\tiny{(LBF)}}$. So the power allocation problem for given beamforming directions  can be written as 
\begin{eqnarray}
\mathbf{\alpha}_{1}^\text{\tiny{(LBF)}},\mathbf{\alpha}_{2}^\text{\tiny{(LBF)}},\cdots,\mathbf{\alpha}_{K}^\text{\tiny{(LBF)}} & = & \max_{\mathbf{\alpha}_{1},\cdots,\mathbf{\alpha}_{K}} \sum_{k=1}^{K}w^k \log_2\Bigl(1+\frac{P\alpha_{k}\mathbf{\bar{g}}_{k}^{\text{\tiny{(LBF)}}{\dagger}}(\mathbf{h}_{k}^{\mbox{\tiny{s}}})^{\dagger}\mathbf{h}_{k}^{\mbox{\tiny{s}}}\mathbf{\bar{g}}_{k}^{\text{\tiny{(LBF)}}}}{\sigma_{k}^{2}}\Bigr),\nonumber \\
\text{subject to: } & \text{ } & \sum_{k=1}^{K}\alpha_{k}\mathbf{\bar{g}}_{k}^{\text{\tiny{(LBF)}}{\dagger}}\mathbf{R}_k^j\mathbf{\bar{g}}_{k}^{\text{\tiny{(LBF)}}}\leq\gamma^j_\text{th},~~ \text{for } j=1,\cdots, J.\label{opt_for_alpha}\end{eqnarray}
The Lagrange function of the above optimization problem   can be written as
\begin{equation}
\mathcal{L}=\sum_{k=1}^{K}w^k \log_2\Bigl(1+\frac{P\alpha_{k}\mathbf{\bar{g}}_{k}^{\text{\tiny{(LBF)}}{\dagger}}(\mathbf{h}_{k}^{\mbox{\tiny{s}}})^{\dagger}\mathbf{h}_{k}^{\mbox{\tiny{s}}}\mathbf{\bar{g}}_{k}^{\text{\tiny{(LBF)}}}}{\sigma_{k}^{2}}\Bigr)-\sum_{j=1}^{J}\biggl(\lambda^j\bigl(\sum_{k=1}^{K}\alpha_{k}\mathbf{\bar{g}}_{k}^{\text{\tiny{(LBF)}}{\dagger}}\mathbf{R}_k^j\mathbf{\bar{g}}_{k}^{\text{\tiny{(LBF)}}}-\gamma^j_\text{th}\bigr)\biggr),\end{equation}
where $\{\lambda^{1}, \cdots, \lambda^{J}\}$ are the Lagrange multipliers. Using KKT conditions, we can write
\begin{equation}
w^k {\Bigl(\alpha_k+\frac{\sigma_{k}^{2}}{P\mathbf{\bar{g}}_{k}^{\text{\tiny{(LBF)}}{\dagger}}(\mathbf{h}_{k}^{\mbox{\tiny{s}}})^{\dagger}\mathbf{h}_{k}^{\mbox{\tiny{s}}}\mathbf{\bar{g}}_{k}^{\text{\tiny{(LBF)}}}}\Bigr)}^{-1}-\sum_{j=1}^{J}\biggl(\lambda^j\mathbf{\bar{g}}_{k}^{\text{\tiny{(LBF)}}{\dagger}}\mathbf{R}_k^j\mathbf{\bar{g}}_{k}^{\text{\tiny{(LBF)}}}\biggr)=0 \;\;\;\;\;\;       \text{for } k=1,\cdots,K, \label{water_filling}
\end{equation}
\begin{equation}
\lambda^j\bigl(\sum_{k=1}^{K}\alpha_{k}\mathbf{\bar{g}}_{k}^{\text{\tiny{(LBF)}}{\dagger}}\mathbf{R}_k^j\mathbf{\bar{g}}_{k}^{\text{\tiny{(LBF)}}}-\gamma^j_\text{th}\bigr)=0,\;\;\;\;\;\;\;\;\;\;\;\;       \text{for } j=1,\cdots,J,\label{slack}\end{equation}
\begin{equation}
\sum_{k=1}^{K}\alpha_{k}\mathbf{\bar{g}}_{k}^{\text{\tiny{(LBF)}}{\dagger}}\mathbf{R}_k^j\mathbf{\bar{g}}_{k}^{\text{\tiny{(LBF)}}}-\gamma^j_\text{th}\leq0,\;\;\;\;\;\;\;\;\;\;\;\;       \text{for } j=1,\cdots,J,\label{constraint_alpha}\end{equation}
\begin{equation}
\lambda^{1}, \cdots, \lambda^{J}\geq0,\end{equation}

 According to eq. (\ref{water_filling}), the power allocation for beamforming direction  corresponding to secondary destination d$_k$ is given by
\begin{equation}
\alpha_{k}^\text{\tiny{(LBF)}}=\max\Biggl(0,\frac{w^k}{\sum_{j=1}^{J}\biggl(\lambda^j\mathbf{\bar{g}}_{k}^{\text{\tiny{(LBF)}}{\dagger}}\mathbf{R}_k^j\mathbf{\bar{g}}_{k}^{\text{\tiny{(LBF)}}}\biggr)}-\frac{\sigma_{k}^{2}}{P\mathbf{\bar{g}}_{k}^{\text{\tiny{(LBF)}}{\dagger}}(\mathbf{h}_{k}^{\mbox{\tiny{s}}})^{\dagger}\mathbf{h}_{k}^{\mbox{\tiny{s}}}\mathbf{\bar{g}}_{k}^{\text{\tiny{(LBF)}}}}\Biggr),\label{capped}
\end{equation}
for $k=1,\cdots,K$. The power allocation in eq. (\ref{capped}) is the cap-limited water-filling solution. 
%
In eq. (\ref{capped}), the power allocation values  $\alpha_{k}^\text{\tiny{(LBF)}}$  are expressed in terms of Lagrange multipliers $\lambda_j$ $(j=1,\cdots, J$) which need to be evaluated.  

In order to obtain the Lagrange multipliers  and consequently $\alpha_{k}^\text{\tiny{(LBF)}}$,  a recursive technique  is used   as  described below. 
First, we assume that only one Lagrange multiplier is greater then zero, i.e., $\lambda^j>0$, while $\lambda^i=0$,  for all $i$ except  $i\neq j$.  This implies that the optimum power allocation values $\alpha_k^\text{\tiny{(LBF)}}$,  $(k=1,\cdots,K)$, satisfy the interference threshold with equality only at  primary receiver p$_j$. For this case, we can write  
\begin{equation}
\sum_{k=1}^{K}\alpha_{k}^\text{\tiny{(LBF)}}\mathbf{\bar{g}}_{k}^{\text{\tiny{(LBF)}}{\dagger}}\mathbf{R}_k^j\mathbf{\bar{g}}_{k}^{\text{\tiny{(LBF)}}}-\gamma^j_\text{th}=0.
\label{eq_with_one_zero}\end{equation}
Now  the value of $\lambda^j$ and the power allocation values $\alpha_k^\text{\tiny{(LBF)}}$, for all $k$  are found by solving set of equations  in (\ref{capped}) and (\ref{eq_with_one_zero}) simultaneously.
If these values of $\alpha_k^\text{\tiny{(LBF)}}$  satisfy the remaining $(J-1)$ interference constraints given by the set of equations  in (\ref{constraint_alpha}), then $\alpha_k^\text{\tiny{(LBF)}}$ for all $k$  represent the  optimum solution of (\ref{opt_for_alpha}).  Otherwise, we set  $\lambda^{k}>0$ ($k\neq j)$   while $\lambda^i=0$, for all $i$  except  $i\neq k$, and so on  until we find the power allocation values that satisfy all constraints simultaneously. 


If no  power allocation values that satisfy all constraints simultaneously is found, considering one constraint as equality  constraint  we consider the case when two constraints are met with equality. In other words, we set  simultaneously  $\lambda^{j} >0$ and $\lambda^{l}  >0$  while $\lambda^i=0$, for all $i$ except $i\neq j,l$. Then, the following two slackness conditions in  eq. (\ref{slack}) are satisfied as follows
\begin{eqnarray}
\sum_{k=1}^{K}\alpha_{k}^\text{\tiny{(LBF)}}\mathbf{\bar{g}}_{k}^{\text{\tiny{(LBF)}}{\dagger}}\mathbf{R}_k^j\mathbf{\bar{g}}_{k}^{\text{\tiny{(LBF)}}}-\gamma^j_\text{th}=0,
\label{eq_with_two_zeros_1}\\
\sum_{k=1}^{K}\alpha_{k}^\text{\tiny{(LBF)}}\mathbf{\bar{g}}_{k}^{\text{\tiny{(LBF)}}{\dagger}}\mathbf{R}_k^l\mathbf{\bar{g}}_{k}^{\text{\tiny{(LBF)}}}-\gamma^l_\text{th}=0.
\label{eq_with_two_zeros_2}\end{eqnarray}
The values of $\lambda^j$, $\lambda^l$ and the power allocation values $\alpha_k^\text{\tiny{(LBF)}}$ for all $k$  are found by solving the set of  equations in  (\ref{capped}), (\ref{eq_with_two_zeros_1}), and (\ref{eq_with_two_zeros_2}) simultaneously.
If these values of $\alpha_k^\text{\tiny{(LBF)}}$  satisfy the remaining $(J-2)$  interference constraints given by the set of equations in  (\ref{constraint_alpha}), then  $\alpha_k^\text{\tiny{(LBF)}}$ for all $k$  are the  optimum power allocation values. Otherwise, we set another set of two constraints as equality constraint, i.e.,  $\lambda^{m}>0$ and $\lambda^{n}>0$ ($m,n\neq j, l)$ while $\lambda^i=0$, for all $i$ except  $i\neq m,n$,  and so on until we find the values of $\alpha_k^\text{\tiny{(LBF)}}$ that satisfy all constraints simultaneously.  The worst case scenario in terms of complexity occurs when the $J$ constraints hold with equality simultaneously.

This procedure is summarized  below:
\begin{algorithmic}
 \For{$i=1\to J$}
\State - Form $\binom{J}{i}$ different sets, such that each set $\mathcal{S}_{k}^{i}\text{ for }k=1,\cdots,\binom{J}{i}$ is composed of $i$ \State   different $\lambda$'s.
\For {$j=1\to \binom{J}{i}$}
\State-Assume that $\lambda^m=0 \text{ for } \lambda^m \not\in \mathcal{S}_j^{i}$, and that  $\lambda^n>0$ for  $\lambda^n \in \mathcal{S}_j^{i}$, which implies that 
\State    the interference constraints  at $i$ primary receivers are satisfied with equality
\State   simultaneously.
\State -Substitute these $\lambda$'s in eq. (\ref{capped}), and in the slackness conditions given in eqs. (\ref{slack}) to
\State   get the optimum power allocation,   $\alpha_{k}^\text{\tiny{(LBF)}}$ for all $k$.
\State -Check whether the total interference introduced due to  the transmissions to  $K$ seco-
\State ndary destinations  satisfies  the other $(J - i)$ interference constraints given in  eqs. (\ref{constraint_alpha}),
\State - if yes, exit. Otherwise, continue.
 \EndFor
 \EndFor
\end{algorithmic} 

\subsection{Low Complexity Power Allocation Scheme}
The complexity of the power allocation scheme proposed in Section \ref{multiple_SDs} can, in the worst case scenario, be in the order of $\frac{J(J+1)}{2}$. 
The optimum power allocation (OPA) scheme, proposed in Section \ref{multiple_SDs}, jointly finds all the $K$ allocated power values which can, in the worst case, require solving the $J$ interference constraints simultaneously.  Therefore, we also propose a low complexity power allocation (LCPA)  scheme  as described below.

Rather than finding the power allocation value $\mathbf{\alpha}_k$ by keeping all the $J$ interference constraints simultaneously in eq. (\ref{opt_for_alpha}), we propose to find the power allocation value for only  one interference constraint e.g., $j$th interference constraint at a time.  For notational convenience let us denote, the corresponding power value by   $\mathbf{\alpha}_k^{j,\text{\tiny{LCPA}}}$ $(k={1,\cdots,K})$ which can be  written as follows
\begin{equation}
\alpha_{k}^{j, \text{\tiny{LCPA}}}=\max\Biggl(0,\frac{w^k}{\lambda^j\mathbf{\bar{g}}_{k}^{\text{\tiny{(LBF)}}{\dagger}}\mathbf{R}_k^j\mathbf{\bar{g}}_{k}^{\text{\tiny{(LBF)}}}}-\frac{\sigma_{k}^{2}}{P\mathbf{\bar{g}}_{k}^{\text{\tiny{(LBF)}}{\dagger}}(\mathbf{h}_{k}^{\mbox{\tiny{s}}})^{\dagger}\mathbf{h}_{k}^{\mbox{\tiny{s}}}\mathbf{\bar{g}}_{k}^{\text{\tiny{(LBF)}}}}\Biggr).
\label{separate_alphas}
\end{equation}
The value of $\lambda^{j}$ is found from the following complementary slackness condition
\begin{equation}
\lambda^j\bigl(\sum_{k=1}^{K}\alpha_{k}^{j, \text{\tiny{LCPA}}}\mathbf{\bar{g}}_{k}^{\text{\tiny{(LBF)}}{\dagger}}\mathbf{R}_k^j\mathbf{\bar{g}}_{k}^{\text{\tiny{(LBF)}}}-\gamma^j_\text{th}\bigr)=0,\;\;\;\;\;\;\;\;\;\;\;\;       \text{for } j=1,\cdots,J.\end{equation}

So, now for a given  beamforming direction correspond a particular secondary destination d$_k$,   we have $J$ power values $\alpha_{k}^{j, \text{\tiny{LCPA}}}$ $(j=1,\cdots,J)$ corresponding to $J$ interference constraints. Out of these $J$ power values, the minimum power value is selected as the final power allocation value for $k$th beamforming direction, i.e., 
\begin{equation}
\mathbf{\alpha}_k^{\text{\tiny{LCPA}}}=\min\bigl(\mathbf{\alpha}_k^{1, \text{\tiny{LCPA}}}, \mathbf{\alpha}_k^{2, \text{\tiny{LCPA}}},\cdots, \mathbf{\alpha}_k^{J, \text{\tiny{LCPA}}} \bigr).
\label{sub_min}
\end{equation}
The complexity of this proposed LCPA  scheme is in the order of $J$,  compared to that of the OPA scheme which is in the order of $\frac{J(J+1)}{2}$ in the worst case.  This lower complexity comes at the expense of  sum transmission rate of secondary destinations.




\section{\label{Sec:imperfect_channel}Beamforming with Partial Channel Knowledge}
 In many scenarios, the instantaneous CSI of the channels between  the CCRNs and the primary receivers may  not be available at the CCRNs.  During the design process of the cooperative transmit beamforming,  we need to account for the effects of  partial channel knowledge at the CCRNs to ensure a robust protection to  the primary receivers. In this section, we consider two scenarios of having partial  channel knowledge.  The first scenario is  having the erroneous CSI of the channels between the primary users and the CCRNs   due to the imperfect channel estimation.  
 The second scenario is  having only the statistical   CSI of the channels between  the primary users and the CCRNs  rather than the instantaneous CSI.  For these scenarios our goal is to design  RLBF  techniques.  
\subsection{Beamforming with Erroneous Channel Estimate}

 When the CCRNs  have erroneous estimation of the channels between the primary receivers and the CCRNs, in order to design the RLBF  technique for such scenario, we adopt the following channel estimation uncertainty model. If the channel estimation of $\mathbf{h}_j^{\mbox{\tiny{P}}}$ is erroneous, the estimation error can be modeled as 
\begin{equation}
\mathbf{h}_{j}^{\mbox{\tiny{p}}}=\hat{\mathbf{h}}_{j}^{\mbox{\tiny{p}}}+\mathbf{e}_j,
\label{error_channel}\end{equation}
where $\mathbf{h}_{j}^{\mbox{\tiny{p}}}$ is the actual instantaneous channel vector between the CCRNs  and   primary receiver p$_j$, $\hat{\mathbf{h}}_{j}^{\mbox{\tiny{p}}}$  is the estimated channel vector between the CCRNs and  primary receiver, p$_j$  and $\mathbf{e}_j$ is the corresponding estimation error vector. Based on the accuracy of the estimation technique used, the channel estimation uncertainty can be modeled by the so-called bounded uncertainty model. The bounded uncertainty  model is a well-accepted model that has been used in \cite{uncertainty_model, uncertainty_sphere, Lampe_paper, ref_9_Lampe}. It considers that the uncertainty in the channel estimation is described by a bounded region whose shape depends on the channel estimation technique used. However, a spherical uncertainty region gives the worst case estimation error model \cite{uncertainty_sphere}. In this case, the estimation error vector is bounded by ${\left\Vert{\mathbf{e}_j}\right\Vert}^{2} \leq \epsilon$. 



Using the error model in eq. (\ref{error_channel}), the covariance matrix corresponding to $\mathbf{h}_{j}^{\mbox{\tiny{p}}}$, $\mathbf{R}^j_k$  can be written as \cite{my_journal}
\begin{equation}
\mathbf{R}^j_{k}=\hat{\mathbf{R}}^j_{k}+\Delta^j_{R},
\end{equation}
where $\hat{\mathbf{R}}^j_{k}$ is the estimated covariance matrix corresponding to the  estimated channel fading gains between the  CCRNs     and    primary receiver p$_j$ and  can be calculated using the estimated CSI  $\hat{\mathbf{h}}_j$ as well as  ${\beta}_k^{j(r,f)}$ (see eq. (\ref{chCovMatrix})).   $\Delta^j_{R}$ is the covariance error matrix 
\begin{equation}
\Delta_{R}^j\negthickspace=\negthickspace
\left[\begin{array}{ccc}
{\beta_k^j}^{(1,1)}(e_{j1})^{\dagger}e_{j1} & \negthickspace\cdots\negthickspace & {\beta_k^j}^{(1,L)}(e_{j1})^{\dagger}e_{jL}\\
\negthickspace\vdots\negthickspace & \negthickspace\ddots\negthickspace & \negthickspace\vdots\negthickspace\\
{\beta_k^j}^{(L,1)}(e_{jL})^{\dagger}e_{j1} & \negthickspace\cdots\negthickspace & {\beta_k^j}^{(L,L)}(e_{jL})^{\dagger}e_{jL}\end{array}\right],
\end{equation}
 and is bounded by $\left\Vert \Delta_{R}^j\right\Vert \leq\Psi^j_{R}$, where $\Psi^j_{R}$ is the bound of the uncertainty region of $\hat{\mathbf{R}}^j_{k}$. Since $\mathbf{R}^j_{k}$ is a covariance matrix, it can be factorized using Cholesky decomposition \cite{cholesky}. Therefore, we can write  $\mathbf{R}^j_{k}=\mathbf{C}^j_{k}(\mathbf{C}^j_{k})^{\dagger}$, where $\mathbf{C}^j_{k}$ is a lower triangular matrix. Similarly, we can write $\hat{\mathbf{R}}^j_{k}=\hat{\mathbf{C}}^j_{k}(\hat{\mathbf{C}}^j_{k})^{\dagger}$. Then, the relation between $\mathbf{C}^j_{k}$ and $\hat{\mathbf{C}}^j_{k}$ can be written as
\begin{equation}
\mathbf{C}^j_{k}=\hat{\mathbf{C}}^j_{k}+\Delta^j_{C},\;\;\;\;\;\;\;\;\;\;\;\;\left\Vert \Delta^j_{C}\right\Vert \leq\Psi^j_{C},\label{relation of C}\end{equation}
where $\Psi^j_{C}$ is the bound of the uncertainty region of $\hat{\mathbf{C}}^j_{k}$.

The total asynchronous interference at primary receiver p$_j$ should satisfy the following condition
\begin{equation}
\sum_{k=1}^{K}\mathbf{g}_k^{\dagger}\mathbf{R}^j_{k}\mathbf{g}_k\leq\gamma^j_\text{th}, \text{ for }j=1,\cdots,J.\label{constraint_error}
\end{equation}
Using eq. (\ref{relation of C}) in eq. (\ref{constraint_error}), we can reformulate the total asynchronous interference constraint for   primary receiver p$_j$    as follows
\begin{equation}
\sum_{k=1}^{K}\left\Vert \mathbf{g}_k^{\dagger}{\mathbf{C}^j_{k}}\right\Vert ^{2}\leq\gamma^j _\text{th}, \text{ for }j=1,\cdots,J.\end{equation}
However, in order to ensure a robust design of the beamforming vector using $\mathbf{C}^j_{k}$, the above constraint must be satisfied for the worst case estimate of $\mathbf{C}^j_{k}$, i.e.,
\begin{equation}
\max_{\left\Vert \Delta^j_{C}\right\Vert}\sum_{k=1}^{K}\left\Vert \mathbf{g}_k^{\dagger}{\mathbf{C}}^j_{k}\right\Vert  \leq\sqrt{\gamma^j _\text{th}}.\label{eq:worst_case}\end{equation}
Using the triangle inequality  
and applying Cauchy-Schwartz inequality \cite{my_journal}, we can write
\begin{equation}
\left\Vert \mathbf{g}_k^{\dagger}\mathbf{C}^j_k\right\Vert \leq\left\Vert \mathbf{g}_k^{\dagger}\hat{\mathbf{C}}^j_{k}\right\Vert +\left\Vert \mathbf{g}_k\right\Vert \left\Vert \Delta_{C}^j\right\Vert.\label{cauchy} \end{equation}
Using the maximum value of $\left\Vert \mathbf{g}_k^{\dagger}\hat{\mathbf{C}}^j_k\right\Vert $
given in eq. (\ref{cauchy}) and substituting it in eq. (\ref{eq:worst_case}), the design constraint now becomes
\begin{equation}
\sum_{k=1}^{K}\left\Vert \mathbf{g}_k^{\dagger}\hat{\mathbf{C}}^j_k\right\Vert ^{2}\leq\biggl(\sqrt{\gamma^j _\text{th}}-\sum_{k=1}^{K}\left\Vert \mathbf{g}_k\right\Vert \Psi^j_{C}\biggr)^{2}.\label{in_terms_of_C}\end{equation}

By using the relation between $\hat{\mathbf{R}}^j_k$ and $\hat{\mathbf{C}}^j_{k}$, the design constraint in eq. (\ref{in_terms_of_C}) can finally be expressed as 
\begin{equation}
\sum_{k=1}^{K}\mathbf{g}_k^{\dagger}\hat{\mathbf{R}}^j_k\mathbf{g}_k\leq\biggl(\sqrt{\gamma^j _\text{th}}-\sum_{k=1}^{K}\left\Vert \mathbf{g}_k\right\Vert \Psi^j_{C}\biggr)^{2}.\end{equation}
Therefore, our primal optimization problem for this RLBF technique can be written as 
\begin{eqnarray}
\hat{\mathbf{g}}_{1},\hat{\mathbf{g}}_{2},\cdots,\hat{\mathbf{g}}_{K} & = &  \max_{{\mathbf{g}}_{1},\cdots,{\mathbf{g}}_{K}} \sum_{k=1}^{K}w^k r_{k},\nonumber \\
\text{subject to:} & \text{ } & \sum_{k=1}^{K}\mathbf{g}_k^{\dagger}\hat{\mathbf{R}}^j_k\mathbf{g}_k\leq I^j_\text{th}, \text{ for }j=1,\cdots,J.\label{robust_g}\end{eqnarray}
where $I^j_\text{th}=\biggl(\sqrt{\gamma^j _\text{th}}-\sum_{k=1}^{K}\left\Vert \mathbf{g}_k\right\Vert \Psi^j_{C}\biggr)^{2}$.  The optimization problem in eq. (\ref{robust_g}) is a non-linear and non-convex problem which cannot be solved  optimally. However,  using a two-step  procedure similar to the one described in Section \ref{coopLBF}, a suboptimal solution for the cooperative RLBF can be obtained for the scenario when CCRNs have imperfect CSI of the primary receivers. 


\subsection{Beamforming with Channel Statistics}
When the CCRNs  have the  statistical CSI\footnote{Statistical CSI refers to distribution of CSI which is assumed to be Rayleigh and corresponding  parameter.}  of the channels between the primary receivers and the CCRNs,  the interference thresholds at the primary receivers can be guaranteed statistically. In absence of instantaneous CSI of the channel between primary receiver and a CR transmitter, such statistical interference constraint to primary receivers has been used in   \cite{kim}, \cite{stat_channel}. According to this statistical asynchronous interference constraint,  interference thresholds   are met probabilistically as follows
\begin{equation}
\mbox{Pr}\left(P_\text{asynch}^j\geq\gamma^j_\text{th}\right)\leq\epsilon^j,\label{eq:prob_int}\end{equation}
where Pr denotes probability and $\epsilon^j$ is the maximum allowable probability of violating the interference threshold $\gamma^j_\text{th}$  at primary receiver p$_j$. Since the distribution  of the random interference power $P_\text{asynch}^j$ is not available in a closed-form, the probability in the left side of eq. (\ref{eq:prob_int}) can not be written in a closed-form in terms of average channel gains between the CCRNs and the primary receivers.  In what follows we develop an upper bound on  this  probability value, i.e.,  $\mbox{Pr}\left(P_\text{asynch}^j\geq\gamma^j_\text{th}\right)$, using the well-known Markov's inequality, in terms of average channel fading power gains between  primary receiver p$_j$  and CCRNs. 

According to the Markov's inequality the probability that a nonnegative random variable $X$ is greater than or equal to some positive constant $a$ is upper bounded by the ratio of expected value of $X$ and $a$ i.e., Pr$(X\geq a)\leq \frac{\mbox{E}(X) }{a}$ \cite{Markov_ineq}. Since the asynchronous interference power $P_\text{asynch}^j$ is a non-negative function of the random variables $h_{jr}^{\mbox{\tiny{P}}},\; r=1,\cdots,L$,  according to Markov's inequality, the probability $\mbox{Pr}\left(P_\text{asynch}^j\geq\gamma^j_\text{th}\right)$ is upper bounded as follows
\begin{equation}
\mbox{Pr}\left(P_\text{asynch}^j\geq\gamma^j_\text{th}\right)\leq\frac{\mbox{E}\left(P_\text{asynch}^j\right)}{\gamma^j_\text{th}}
\label{Markov}
\end{equation}
which leads to a limit on the average asynchronous interference power on primary receiver p$_j$  (c.f.  eq. (\ref{eq:prob_int}))
\begin{equation}
\mbox{E}\left(P_\text{asynch}^j\right)\leq\epsilon^j\gamma^j_\text{th}.
\label{ubintcons}
\end{equation}
Since the total asynchronous interference power at primary receiver p$_j$, $P_\text{asynch}^j$, is the summation of the interference powers corresponding to the transmissions of different secondary destinations,  the average value of the total asynchronous interference power at primary receiver p$_j$  can be written as
\begin{equation}
\mbox{E}\left(P_\text{asynch}^j\right)=\sum_{k=1}^{K}\mbox{E}\left(P_\text{asynch}^{(j,k)}\right).\end{equation}

The  interference power  at p$_j$  resulting from transmission to secondary destination d$_k$, $P_\text{asynch}^{(j,k)}$ can be written in expanded form as follows (c.f. eq. (\ref{AsynPower}))
\begin{equation}
P_\text{asynch}^{(j,k)}=\sum_{r=1}^{L}\sum_{\begin{array}{c}
f=1,f\neq r\end{array}}^{L}g_{kf}^{\dagger}(h_{jf}^{\mbox{\tiny{p}}})^{\dagger}h_{jr}^{\mbox{\tiny{p}}} g_{kr}{\beta_k^{j}}^{(r,f)}+\sum_{r=1}^{L}g_{kr}^{\dagger}\left|h^{\mbox{\tiny{p}}}_{jr}\right|^{2}g_{kr}{\beta_k^{j}}^{(r,r)}.\label{signal_i}\end{equation}
Since the channel fading coefficients between different CCRNs and  primary receiver p$_j$  are independent and have zero mean, the average value of the first term in eq. (\ref{signal_i}) is equal to zero. For the second term in eq. (\ref{signal_i}), it can be easily shown that for a Rayleigh fading channel, the   the fading power gain, $\left|h^{\mbox{\tiny{p}}}_{jr}\right|^{2}$ has an exponential distribution with a mean value of $\Omega^j_r$, where $\Omega^j_r$. The term $\sum_{r=1}^{L}g_{kr}^{\dagger}\left|h^{\mbox{\tiny{p}}}_{jr}\right|^{2}g_{kr}{\beta^{j}}^{(r,r)}$  is a summation of $L$ independent and identically distributed (i.i.d.) exponential random variables, which is a hypo-exponential random variable, with a mean value of $\sum_{r=1}^{L}g_{kr}^{\dagger}\Omega^j_{r}g_{kr}{\beta^{j}}^{(r,r)}$.
Therefore the average value of $P_\text{asynch}^{(j,k)}$ is given by
\begin{equation}
\mbox{E}\left(P_\text{asynch}^{(j,k)}\right)=\sum_{r=1}^{L}g_{kr}^{\dagger}\Omega^j_{r}g_{kr}{\beta_k^{j}}^{(r,r)}.\end{equation}
This average interference power at primary receiver  p$_j$ can be rewritten in a matrix form as follows
\begin{equation}
\mbox{E}\left(P_\text{asynch}^{(j,k)}\right)=\mathbf{g}_k^{\dagger}\bar{\mathbf{R}}_k^j\mathbf{g}_k
\label{matrixformint}
\end{equation}
where
\begin{equation}
\bar{\mathbf{R}}_k^j=\left[\begin{array}{ccc}
{\beta_k^{j}}^{(1,1)}\Omega^j_{1} & \cdots & 0\\
\vdots & \ddots & \vdots\\
0  & \cdots & {\beta_k^{j}}^{(L,L)}\Omega^j_{L}\end{array}\right].
\end{equation}
Using eq. (\ref{matrixformint}),  eq. (\ref{ubintcons}) can be written as
\begin{equation}
\sum_{k=1}^{K}\mathbf{g}_k^{\dagger}\bar{\mathbf{R}}_k^j\mathbf{g}_k\leq\epsilon^j\gamma^j_\text{th}.
\label{eq:avg:int}
\end{equation}
Now the cooperative RLBF  vector that maximizes the weighted sum rate of secondary destinations while satisfying the new interference constraint in eq. (\ref{eq:avg:int}) can be formulated as an optimization problem as follows 
%
\begin{eqnarray}
\hat{\mathbf{g}}_{1},\hat{\mathbf{g}}_{2},\cdots,\hat{\mathbf{g}}_{K} & = &  \max_{{\mathbf{g}}_{1},\cdots,{\mathbf{g}}_{K}} \sum_{k=1}^{K}w^k r_{k},\nonumber \\
\text{subject to:} & \text{ } & \sum_{k=1}^{K}\mathbf{g}_k^{\dagger}\bar{\mathbf{R}}_k^j\mathbf{g}_k\leq\epsilon^j\gamma^j_\text{th}, \text{ for }j=1,\cdots,J.\end{eqnarray}

The above optimization problem is again a non-linear and non-convex optimization problem, which cannot be solved optimally. However, using the two steps procedure described in Section \ref{coopLBF}, a suboptimal solution for the cooperative RLBF for the statistical CSI scenario can be obtained by substituting for $\gamma^j_\text{th}$ by $\epsilon^j\gamma^j_\text{th}$ in eq (\ref{optimum_form}).

\section{\label{sec:CCRN-Selection} Joint CCRN Selection and Cooperative Beamforming}

Since different CCRNs are located in different geographical locations, their contributions vary significantly towards the interfering signals at the primary receivers, as well as the received signals at the secondary destinations. Intuitively, a CCRN selection strategy can further improve the performance of the cooperative beamforming algorithm. The joint design of cooperative beamforming and relay selection has been studied before for conventional cooperative networks, see for example \cite{relay_selection} and the references therein. In  our previous work in \cite{my_journal}, we studied relay selection  strategies for  CR systems with a single primary receiver  and a single secondary destination.  In this section, for completeness of the considered system with multiple primary and multiple secondary receivers,  we extend the  joint design of cooperative beamforming and CCRN selection. 

To formulate the CCRN selection problem mathematically, we define a CCRN selection vector $\mathbf{\mathbf{S}}$ of size $K\times 1$, where $K$ is the number of CCRNs in the network. The elements of $\mathbf{S}$,  $s_j$ can take value of either $1$ or $0$, to indicate whether the  CCRN c$_j$ has been selected for transmission or not respectively. For notational convenience,  we define $\mathbf{W}$ as a diagonal matrix having its diagonal elements equal to those of vector $\mathbf{S}$, as follows
\begin{equation}
\mathbf{W}= \text{Diag}(\mathbf{S}).
\end{equation}

 The received signal power at secondary destination d$_k$,  $P_{k,\text{sig}}$ is given by \cite{my_journal} \begin{equation}P_{k,\text{signal}}=P{\mathbf{g}_k}^{\dagger}\mathbf{W}^{\dagger}(\mathbf{h}_k^{\mbox{\tiny{s}}})^{\dagger}\mathbf{h}_k^{\mbox{\tiny{s}}}\mathbf{W}\mathbf{g}_k.
\label{RxPwrRelSel}
\end{equation}
The problem of joint CCRN selection and cooperative beamforming is considered as a  mixed-integer non linear problem (MINLP), since the elements of $\mathbf{W}$ can only take a value of $0$ or $1$. Such MINLP can be solved by decoupling it into a non-linear problem (NLP) and a mixed-integer linear problem (MILP), as in \cite{ref_10_MINLP,ref_13_MINLP}.  

Using the value of the received signal power in eq. (\ref{RxPwrRelSel}) for a given relay selection matrix $\mathbf{W}$, the  beamforming vector $\mathbf{g}_k^\text{\tiny{(S)}}{\scriptstyle{(\mathbf{W})}}$ can be found from the following NLP problem:
\begin{eqnarray}
\mathbf{g}_{1}^\text{\tiny{(S)}}{\scriptstyle{(\mathbf{W})}},\mathbf{g}_{2}^\text{\tiny{(S)}}{\scriptstyle{(\mathbf{W})}},\cdots,\mathbf{g}_{K}^\text{\tiny{(S)}}{\scriptstyle{(\mathbf{W})}}&=&\max_{\mathbf{g}_{1},\cdots,\mathbf{g}_{K}}\sum_{k=1}^{K}w^k \log_2\Biggl(1+\frac{P\mathbf{g}_{k}^{\dagger}\mathbf{W}^{\dagger}(\mathbf{h}_{k}^{\mbox{\tiny{s}}})^{\dagger}\mathbf{h}_{k}^{\mbox{\tiny{s}}}\mathbf{W}\mathbf{g}_{k}}{\sigma_{k}^{2}+{\displaystyle \sum_{
i=1, i\neq k}^{K}\mathbf{g}_{i}^{\dagger}\mathbf{W}^{\dagger}\mathbf{T}^k_{i}\mathbf{W}\mathbf{g}_{i}}}\Biggr),\nonumber \\
\text{subject to: } & \text{ } & \sum_{k=1}^{K}\mathbf{g}_{k}^{\dagger}\mathbf{W}\mathbf{R}_k^j\mathbf{W}\mathbf{g}_{k}\leq\gamma^j_\text{th}, ~~~\text{for } j=1,\cdots, J.\label{eq:constriant_with_Selection}\end{eqnarray}
For a given CCRN selection matrix $\mathbf{W}$,  the optimization problem in eq. (\ref{eq:constriant_with_Selection}) is   non-convex and non-linear similar to the one in eq. (\ref{optimum_form}).  Therefore, for a given CCRN selection matrix $\mathbf{W}$,   we can use the same two-phase suboptimal approach described in Section \ref{coopLBF} to find the suboptimal   ${\mathbf{g}_k}^\text{\tiny{(S,sub)}}{\scriptstyle{(\mathbf{W})}}$ for all $k$. Then, 
the optimal CCRN  selection matrix $\mathbf{W^{\ast}}$  and  corresponding  ${\mathbf{g}_k}^\text{\tiny{(S,sub)}}(\scriptstyle{\mathbf{W^{\ast}}})$ are obtained via exhaustive search over all possible selections of  $\mathbf{W}$.



The CCRN selection scheme can be further extended to  the case where only partial knowledge of the channels between the primary receivers and the CCRNs is available at the CCRNs. Using similar procedures as those described in Section \ref{Sec:imperfect_channel}, we can jointly design the CCRN selection and the cooperative beamforming for the case of erroneous channel estimate, and having only the statistical CSI available at the CCRNs. Due to space limitation, we do not include them. 

\vspace{-0.1in}
\section{\label{sec:Numerical-results}Numerical results }
\begin{figure}[h]
\centering
\includegraphics[scale=0.5,angle=-90,trim=2cm 0cm 2.5cm 0cm]{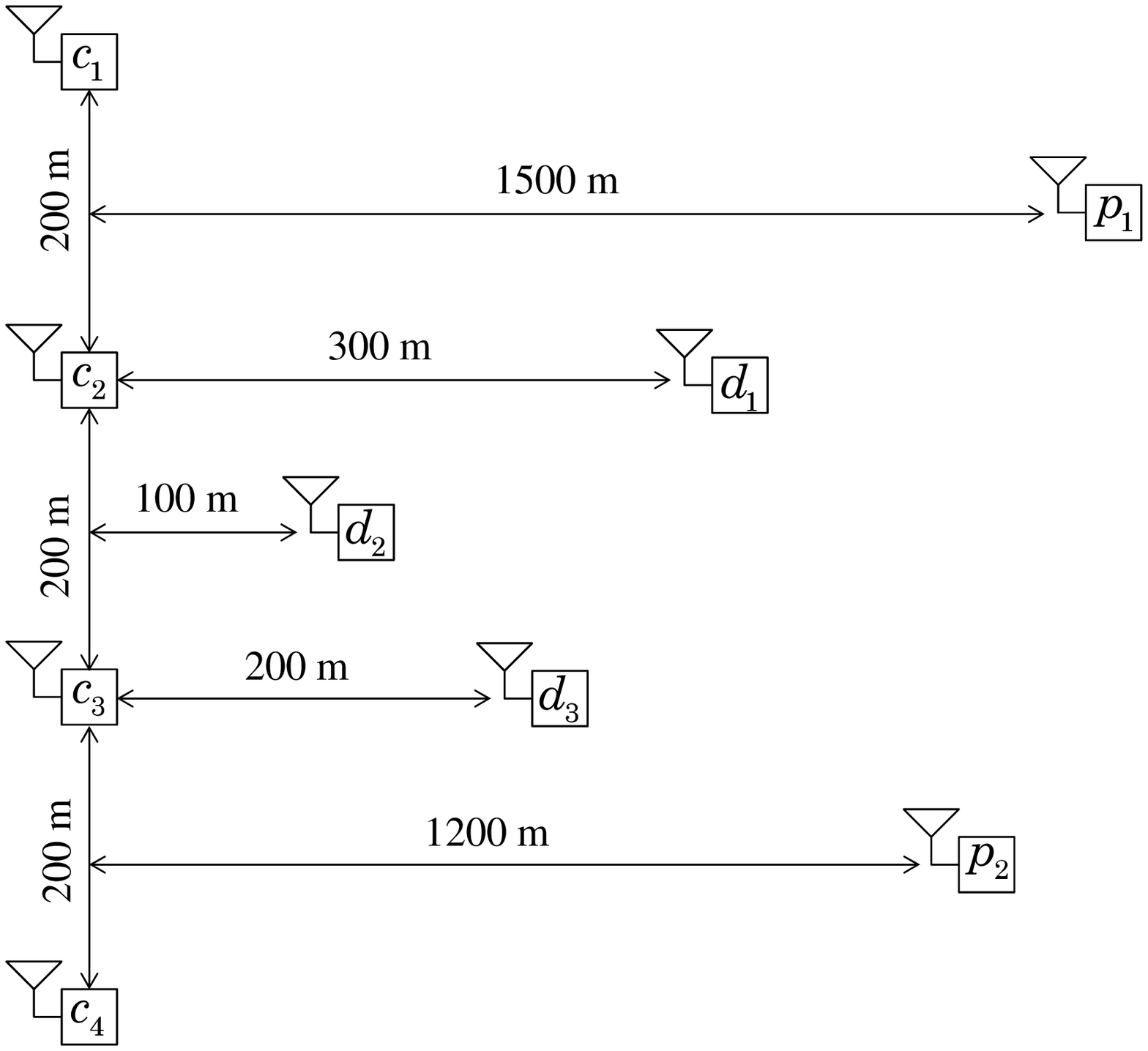}
\caption{Simulated network topolgy.}
\label{network_model}
\end{figure}

In this section, we present some selected numerical results in order to compare the performances of various beamforming techniques in the presence of asynchronous interference. For all the numerical examples presented in this section, we consider the network topology shown in Fig. \ref{network_model}. For the sake of simplicity, it is assumed to have $J=2$  primary receivers, $K=3$ secondary destinations  and $L=4$ CCRNs.   The distances between the nodes in Fig. \ref{network_model} are picked up arbitrarily. Distances between other nodes and corresponding propagation delays can be obtained easily using the given distances.   We assume that all the channel fading gains are identically and independently Rayleigh distributed. We consider a slot duration $T_\text{slot}=0.4$ msec, and a log-distance path loss model with a path loss exponent value of $4$.  The normalized average interference power from the primary transmitter to the secondary destinations, d$_1$, d$_2$ and d$_3$ are  assumed to be  $-10, -20$, and $-15$ dB,  respectively. For simplicity, we consider   weighting factors  $w_1$, $w_2$ and $w_3$ are equal to one.

\begin{figure}[h]
\centering
\includegraphics[width=0.7\linewidth,height=0.35\textheight, draft=false]{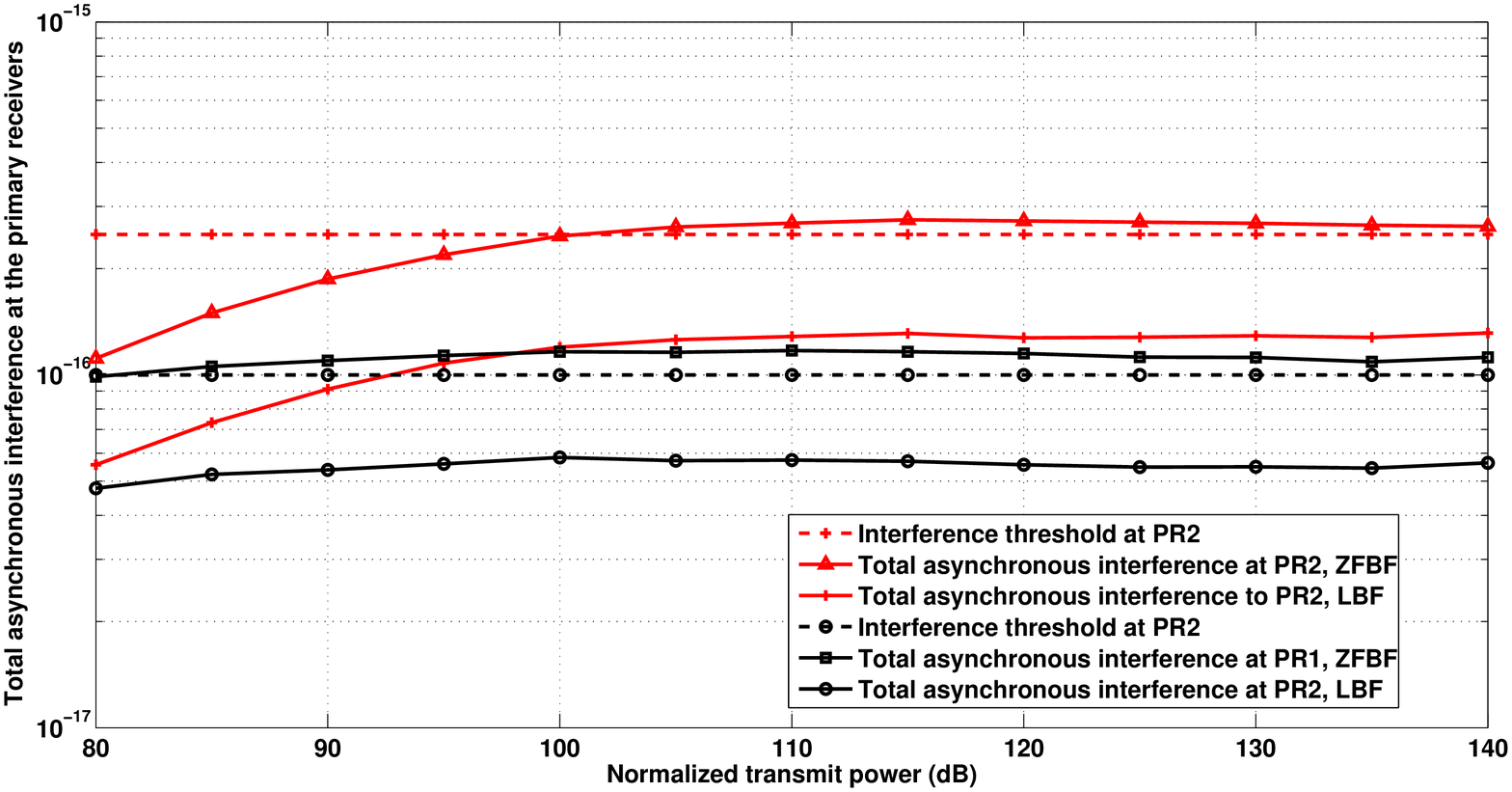}
\vspace{-0.1in}
\caption{Total asynchronous interference power  at the  primary receivers with different interference thresholds ($\gamma^{1}_\text{th}=0.1\times10^{-15}$ and $\gamma^{2}_\text{th}=0.25\times10^{-15}$).
\label{fig:two_users_different_thresholds}}
\end{figure}

In Fig. \ref{fig:two_users_different_thresholds}, we  plot the normalized symbol  power versus the total asynchronous interference signal power introduced at the primary receivers using our proposed cooperative LBF technique. We assume that interference thresholds at primary receiver p$_1$  and p$_2$ are respectively,  $\gamma^{1}_\text{th}=0.1\times10^{-15}$ and $\gamma^{2}_\text{th}=0.25\times10^{-15}$  which are in the order of the noise power.  In this figure we also plot the   asynchronous interference signal powers introduced at the primary receivers when  ZFBF technique  \cite{distributed_beamforming} is used.   This figure clearly shows that our proposed LBF technique can maintain the asynchronous interference thresholds at the primary receivers simultaneously.  On the contrary, the interference caused by the ZFBF technique  exceeds the interference target thresholds at the primary receivers.  This is expected as ZFBF does not take asynchronous interferences into account in its design.

\begin{figure}[h]
\centering
\centering\includegraphics[width=0.7\linewidth,height=0.35\textheight, draft=false]{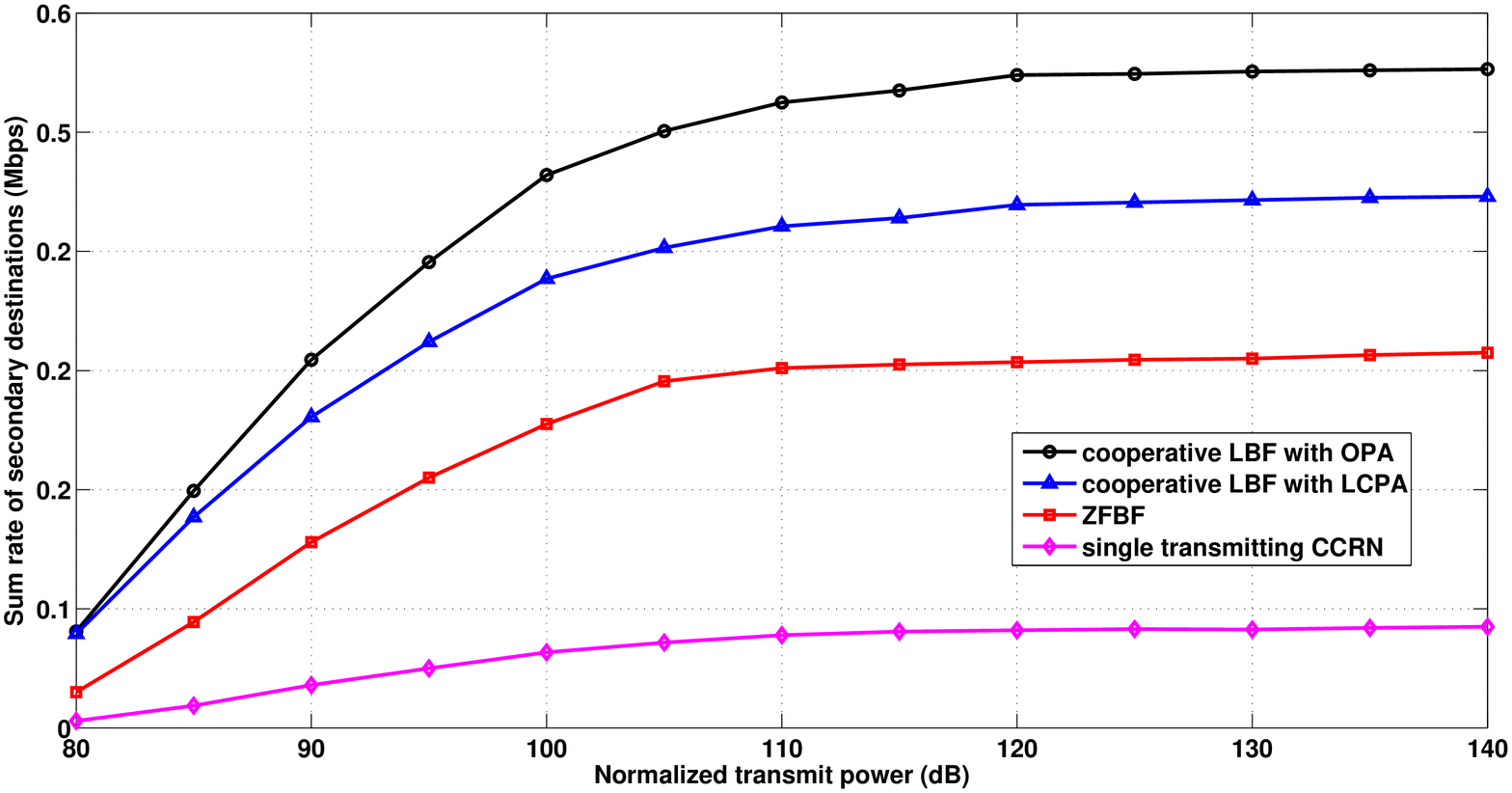}
\vspace{-0.1in}
\caption{Achievable sum transmission rate with various beamforming techniques and single CCRN-based transmission.}
\label{all_schemes}
\end{figure}

In Fig. \ref{all_schemes}, we plot the achievable average sum rate of secondary destinations with the proposed cooperative LBF with OPA scheme,  the LBF with  LCPA  scheme,   and the  ZFBF technique. For the sake of completeness in Fig. \ref{all_schemes}, we also plot the achievable sum  transmission rate of secondary destinations when a single CCRN is selected  for transmission. In this case no beamforming is applied and  we select the CCRN that  maximizes the sum rate of all the secondary destinations. The selected CCRN uses  a transmit power value that satisfies all the primary interference constraints.  With ZFBF an outage is considered    if the instantaneous interference caused by the CCRNs  at any primary receiver exceeds its corresponding target threshold $\gamma^j_\text{th}$. From this figure we can observe that the proposed LBF can achieve a higher sum rate than the well-known ZFBF technique  for the CR-based broadcasting system.   In particular, the increase in sum transmission rate of secondary destinations is about 64\%. This reason can be explained as follows. 
The ZFBF technique  can not satisfy  interference threshold(s)  often and it  leads to a  frequent transmission outage. As such the overall transmission rate of the secondary destinations is degraded.   From Fig. \ref{all_schemes}, we can also see that the proposed  LCPA scheme that  has a lower complexity  suffers from a performance degradation compared to  the OPA scheme as expected.  However, the LBF technique with LCPA scheme  achieves a higher transmission rate compared to the ZFBF technique.  We can also observe from this figure that the single CCRN-based transmission offers the lowest possible transmission rate for the CR system.   This can be explained by the fact that the  single CCRN-based transmission scheme does not take advantage  of beamforming which improves received signal power at the secondary destinations while minimizing the effect of asynchronous interference at the primary receivers.
 
\vspace{-0.1in}

\begin{figure}[h]
\centering
\includegraphics[width=0.7\linewidth,height=0.35\textheight, draft=false]{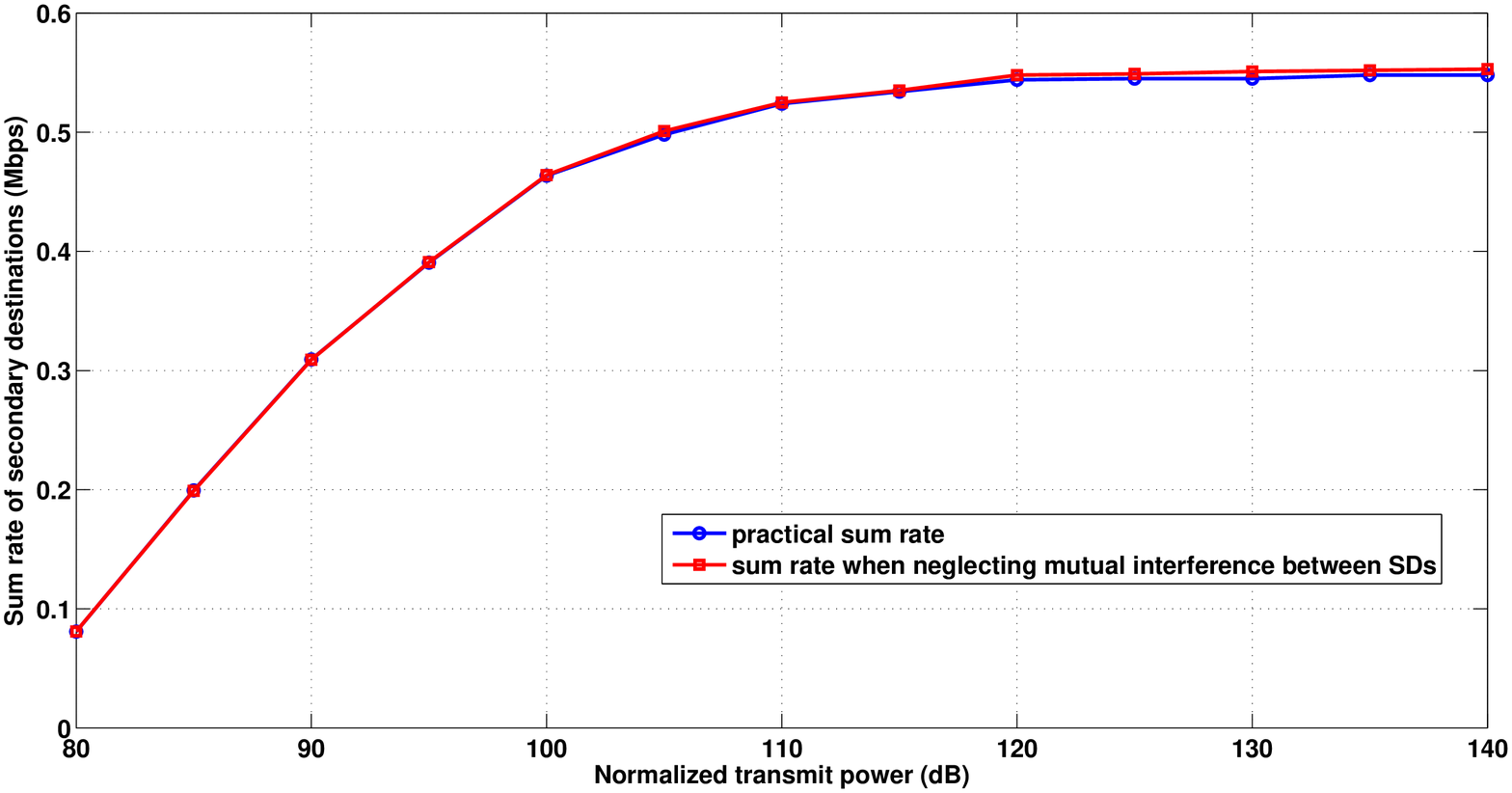}
\vspace{-0.1in}
\caption{Approximated  sum transmission rate versus actual sum transmission rate.}
\label{validate_rate}
\end{figure}

As we have mentioned in Section \ref{coopLBF}  that  after minimizing the asynchronous interference powers,    the mutual asynchronous interferences between secondary destinations can be neglected. Based on this assumption, in Phase-II (see Section III-B2) we have developed the OPA  scheme (or the LCPA scheme) among different beamforming directions.  In order to study the validity   of  such assumption,  in Fig. \ref{validate_rate}   we compare the sum rate of  secondary destinations for two cases. The first case is the practical case in which we calculate the actual sum rate of the secondary destinations taking into account the mutual asynchronous interference signals. In the second case, we use the approximation   (c.f. eq. (\ref{rApprox}))  in which the mutual asynchronous interference signals at the secondary destinations are neglected compared to the noise power.  Fig. \ref{validate_rate}  shows that the approximated and practical values of the sum rate are almost equal. This validates the assumption  of neglecting the mutual asynchronous interferences between secondary destinations in our development of the suboptimal  LBF technique. 

\vspace{-0.1in}

\begin{figure}[h]
\centering
\includegraphics[width=0.7\linewidth,height=0.35\textheight, draft=false]{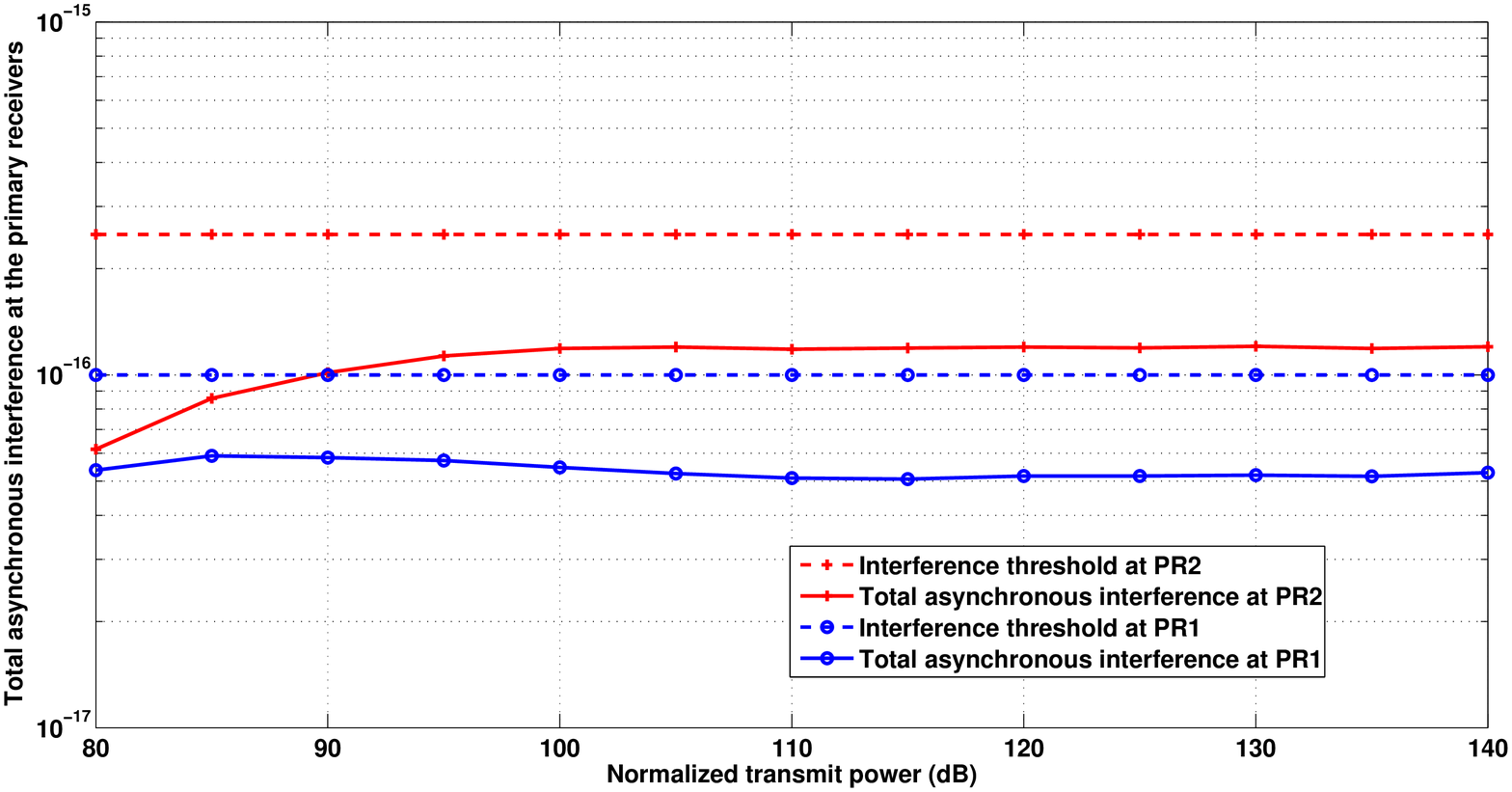}
\caption{Transmit power versus asynchronous  interference powers  at the  primary receivers using the RLBF  technique.}
\label{robust_technique}
\end{figure}

\begin{figure}[h]
\centering
\includegraphics[width=0.7\linewidth,height=0.35\textheight, draft=false]{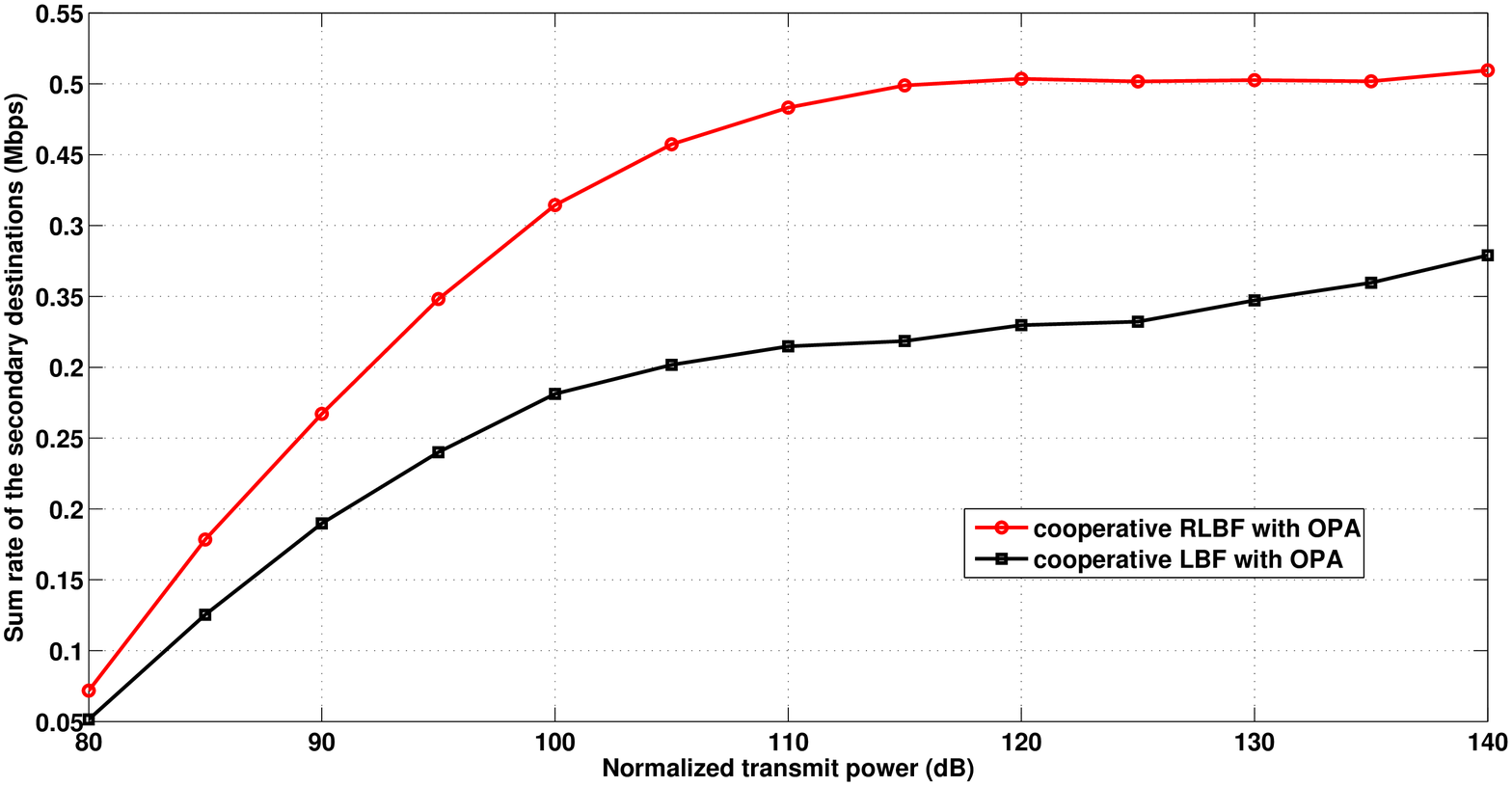}
\caption{The sum transmission  rate of  the  secondary destinations  with LBF and RLBF techniques when there is a certain  channel estimation error.}
\label{robust_unsatisfied}
\end{figure}

Next, we investigate the performances of  the proposed RLBF  techniques in case of having partial CSI  between the primary receivers and CCRNs.  In Fig. \ref{robust_technique}, we plot the   asynchronous interference powers at both primary receivers assuming the erroneous channel estimates at the CCRNs.  From this figure, it is obvious that the RLBF  technique can meet the interference thresholds of the primary receivers even when an erroneous estimation of the channels are available at the CCRNs. We also compare the performance of the LBF technique that requires perfect channel knowledge  with that of the RLBF  technique in Fig. \ref{robust_unsatisfied}. In this figure, we plot the achievable average sum rate of both techniques.  For a fair comparison, with LBF technique, that neglects the channel estimation error at the CCRNs, an outage is considered    if the instantaneous interference caused by the CCRNs  at any primary receiver exceeds its corresponding target threshold $\gamma^j_\text{th}$. We can see from Fig.~\ref{robust_unsatisfied}  that the RLBF  technique achieves a higher sum rate  for secondary destinations  compared to the cooperative LBF technique  when there is a certain channel estimation error. This is expected as the cooperative LBF does not take channel estimation error into account in its design. As such the interference thresholds at the primary receiver(s) can exceed frequently. The value of the estimation error bound,  $\Psi_{C}$  is $0.25\times10^{-8}$. This value has been chosen for the estimation error bound, because lower values will not capture the violation of  the interference threshold with the LBF technique.

\begin{figure}[h]
\centering
\includegraphics[width=0.7\linewidth,height=0.35\textheight, draft=false]{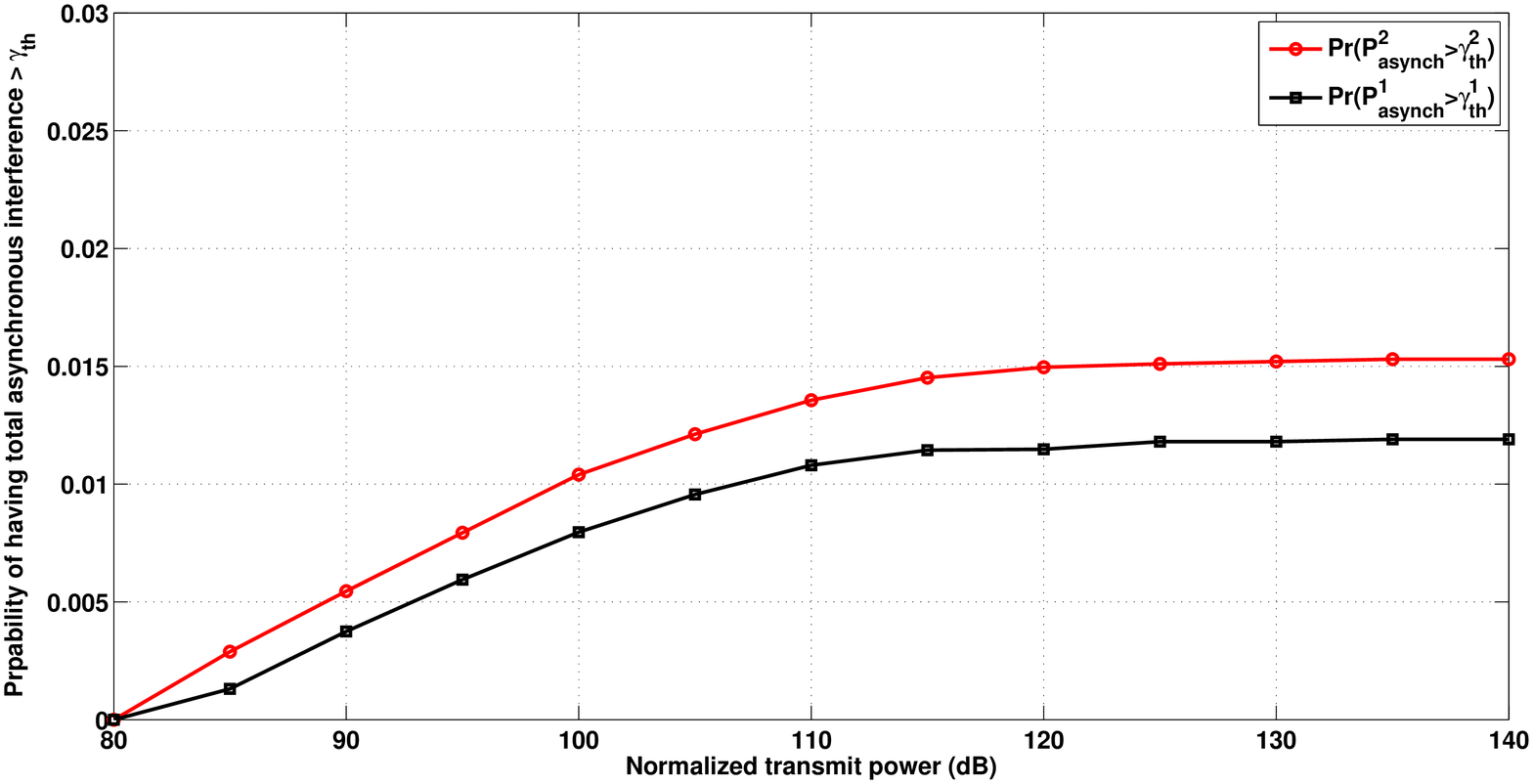}
\caption{The probability that the total asynchronous interference at each primary receiver is greater than $\gamma^j_\text{th}$ when having statistical CSI at the CCRNs.}
\label{violation_prob}
\end{figure}

\begin{figure}[h]
\centering\includegraphics[width=0.7\linewidth,height=0.35\textheight, draft=false]{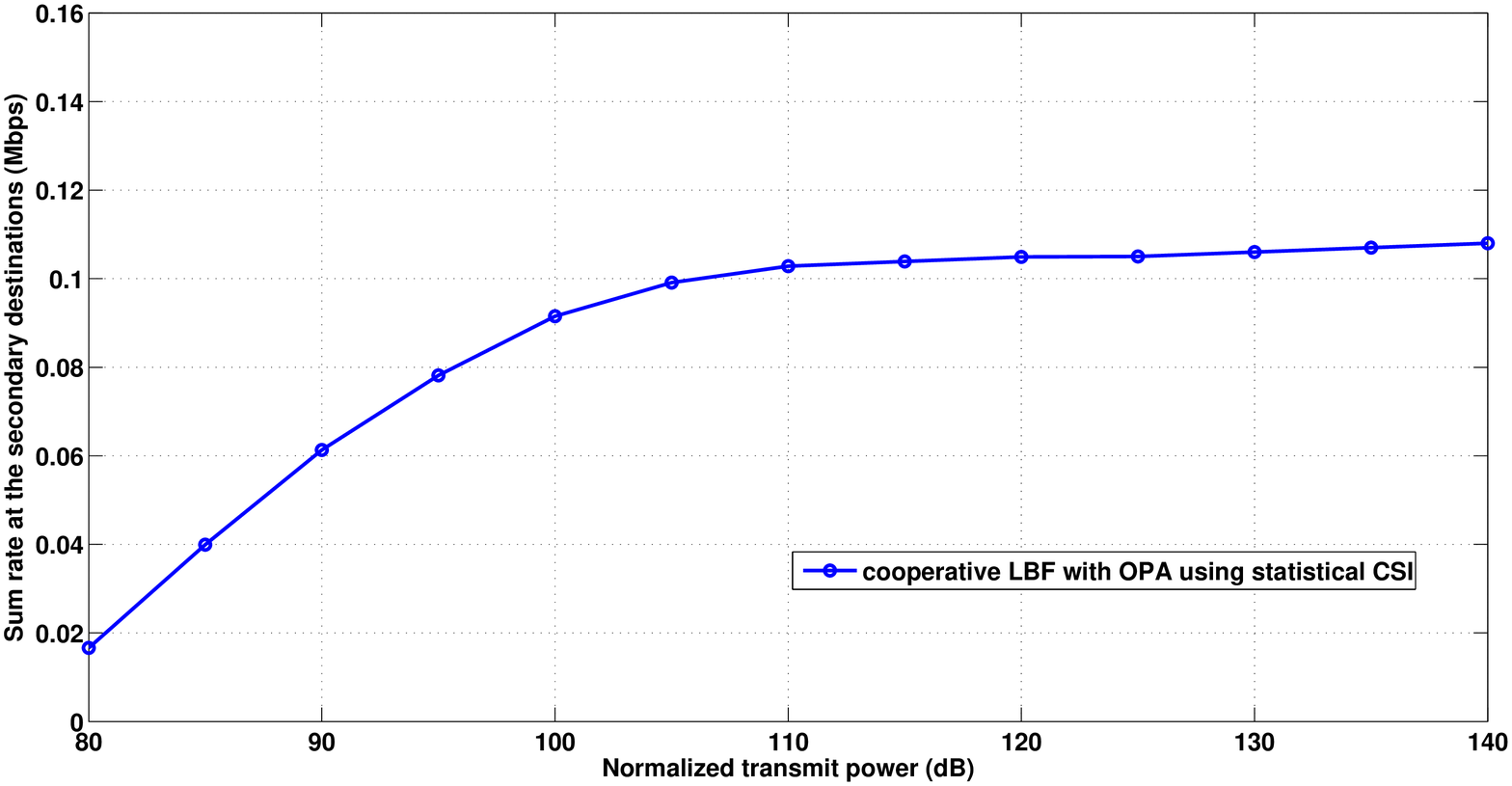}
\caption{The sum rate of  the secondary destinations  with average CSI at CCRNs  and statistical interference constraints at primary receivers.}
\label{statistical}
\end{figure}

We also investigate the performance of  the proposed RLBF technique in case of having only statistical CSI of the  channels between the primary receivers and the CCRNs.  In Fig. \ref{violation_prob}, we plot  the probability of having the instantaneous asynchronous interference power at each primary receiver larger than its target threshold. The value of the maximum allowable probability of violating the interference thresholds $\epsilon^1,\epsilon^2$ are assumed to be $0.1$. It is obvious  from Fig. \ref{violation_prob} that the probability of violating the interference thresholds  is  maintained within the maximum allowable probability value. 
 The achievable average sum rate of this RLBF  technique is shown in Fig. \ref{statistical}. 

\begin{figure}[h]
\centering
\includegraphics[width=0.7\linewidth,height=0.35\textheight, draft=false]{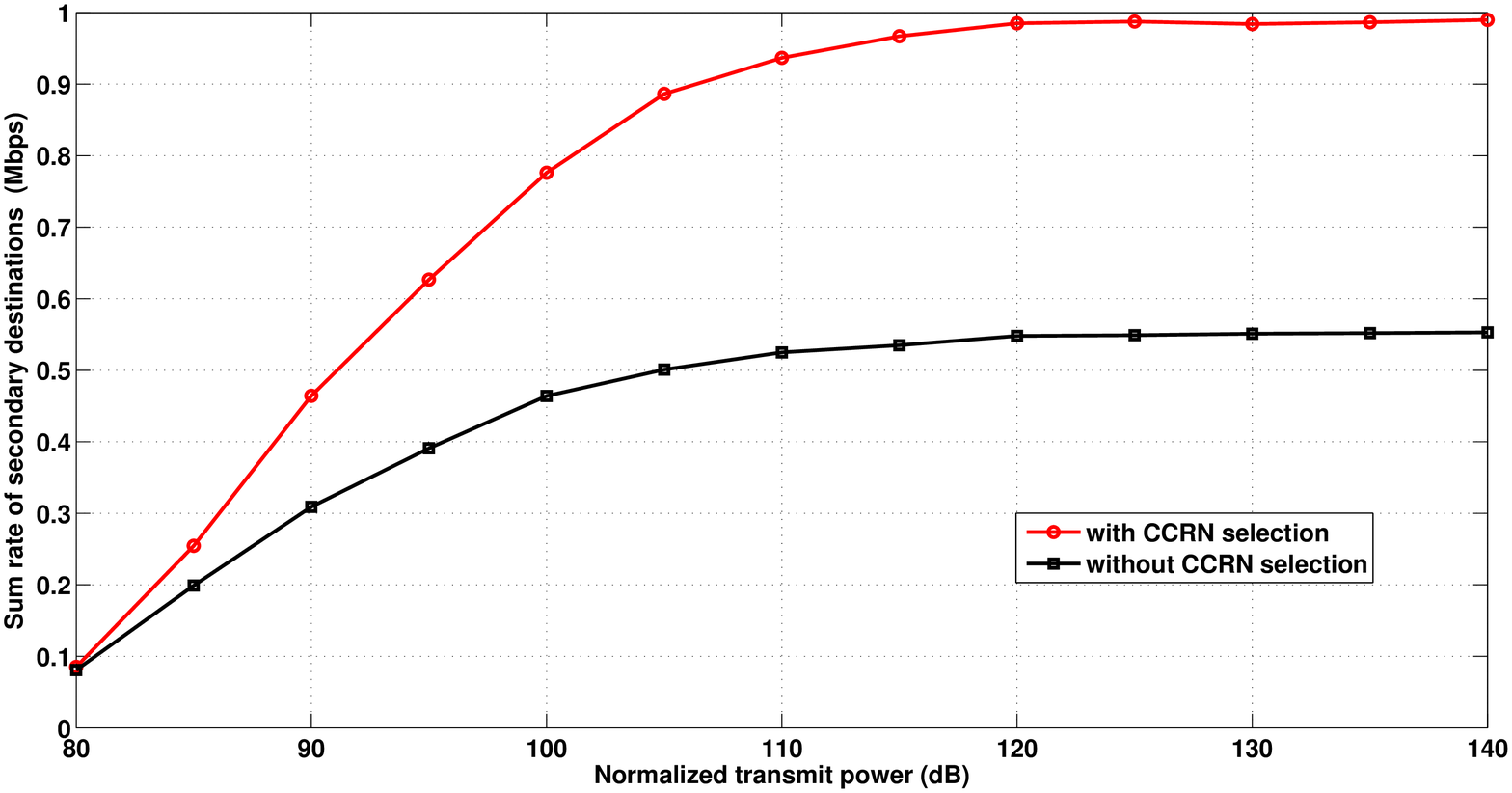}
\caption{The sum rate of  the  secondary destinations  with  and without CCRN selection strategy.}
\label{relay_selection}
\end{figure}

Finally the performance enhancement achieved by applying  the CCRN selection scheme  in conjunction with the LBF technique proposed in section \ref{sec:CCRN-Selection}, is investigated  in Fig. \ref{relay_selection}. In particular, in this figure we plot the average achievable sum rate of the secondary destinations with  joint CCRN selection and cooperative beamforming technique. In this figure we also plot the achievable sum rate of the LBF technique assuming that  all the CCRNs participate in beamforming (i.e.,  without applying any CCRN selection strategy).
From this figure, it is interesting to see that  the CCRN selection scheme in conjunction with the LBF technique outperforms the LBF technique when no CCRN selection is employed.  This increase is about 45\% and  the reason  can be explained intuitively as follows.   When a CCRN selection scheme is employed, the CCRNs are selected judiciously considering their contributions towards the achievable sum rate at the secondary destinations as well as the total interference power at the primary receivers.

\section{\label{sec:conc}Conclusion}

In order to address the  asynchronous interference issue, in this paper, we  have proposed  innovative cooperative beamforming techniques for a generalized CR radio-based broadcasting  system with multiple primary and multiple secondary receivers. In particular,  the cooperative beamforming design is formulated as an optimization problem that maximizes the weighted sum achievable transmission rate  of secondary destinations while it  maintains the asynchronous interferences at the primary receivers below their target thresholds. In light of the intractability of the
problem,  
we have proposed a two-phase suboptimal  beamforming technique.  We have considered both perfect and imperfect CSI of channels between CCRNs and primary receivers.    We also have investigated  the performance of joint CCRN selection and  beamforming technique.  The presented numerical results have shown  that the proposed beamforming technique can significantly reduce the interference signals at all primary receivers and  can provide an increase up to 64\%  in  the sum transmission rate of secondary destinations compared to  the well known  zero-forcing beamforming (ZFBF) technique.  The presented results  have also shown  that   cooperating beamforming node  selection in conjunction with  beamforming can further increase (up to 45\%)  sum data rate of secondary destinations. 
The presented  numerical results have shown that, our proposed robust design of the beamforming vector can  maintain the asynchronous interference constraints at multiple primary receivers when partial CSI is available at the CCRNs.

%
%
\bibliographystyle{ieeetr}
\bibliography{asynch_bib}

\end{document}